\documentclass[]{amsart}
\usepackage{amsmath,amssymb, amscd}
\usepackage{graphicx}
\usepackage{color}

\usepackage{cite}

\begin{document}
\title[Numerical study of DS I]{Numerical study of Davey-Stewartson I 
systems}

\author[J.~Frauendiener]{Joerg Frauendiener}
\address[J.~Frauendiener]{Department of Mathematics and Statistics, 
University of Otago,      
P.O. Box 56, Dunedin 9010, New Zealand}
\email{joergf@maths.otago.ac.nz}

\author[C.~Klein]{Christian Klein}
\address[C.~Klein]{Institut de Math\'ematiques de Bourgogne
9 avenue Alain Savary, BP 47870, 21078 Dijon Cedex}
\email{christian.klein@u-bourgogne.fr}

\author[U.~Muhammad]{Umar Muhammad}
\address[U.~Muhammad]{Institut de Math\'ematiques de Bourgogne, UMR 
5584\\
                Universit\'e de Bourgogne-Franche-Comt\'e, 9 avenue Alain Savary, 21078 Dijon
                Cedex, France}
    \email{umarmuhddauda@gmail.com}

\author[N.~Stoilov]{Nikola Stoilov}

\address[N.~Stoilov]{Institut de Math\'ematiques de Bourgogne, UMR 
5584\\
                Universit\'e de Bourgogne-Franche-Comt\'e, 9 avenue Alain Savary, 21078 Dijon
                Cedex, France}
    \email{Nikola.Stoilov@u-bourgogne.fr}

\begin{abstract}
	An efficient high precision hybrid numerical approach for 
	integrable Davey-Stewartson (DS) I equations for trivial boundary 
conditions at infinity is presented for 
	Schwartz class initial data. 
The code is used for a detailed numerical study of DS I solutions in this 
class. Localized stationary solutions are constructed and 
shown to be unstable against dispersion and blow-up. A finite-time 
blow-up of  initial data in the Schwartz class of smooth rapidly 
decreasing functions is 
discussed.  

\end{abstract}

\date{\today}

\subjclass[2000]{}

\thanks{This work was partially supported  by 
the ANR-FWF project ANuI - ANR-17-CE40-0035, the isite BFC project 
NAANoD, the ANR-17-EURE-0002 EIPHI and by the 
European Union Horizon 2020 research and innovation program under the 
Marie Sklodowska-Curie RISE 2017 grant agreement no. 778010 IPaDEGAN. 
We thanks J.-C.~Saut for helpful discussions and hints}
\maketitle

\section{Introduction} 

This paper is concerned with the numerical study of the  integrable 
Davey-Stewartson (DS) I system  written in the 
form
\begin{equation}
\label{DSIsys}
\begin{array}{ccc}
i
\Psi_{t}+ \Psi_{xx}+\Psi_{yy}+2\left(\Phi+\left| \Psi \right|^{2}\right)\Psi & = & 0,
\\
\Phi_{xx}-\Phi_{yy}+2\left| \Psi \right|_{xx}^{2} & = & 0,
\end{array}
\end{equation}
where  indices denote partial derivatives,  and where $\Phi$ denotes 
a mean field. In the classification of Ghidaglia and Saut \cite{GS}, 
this is an elliptic-hyperbolic equation since the second order 
operator acting on $\Psi$ in 
the first equation of (\ref{DSIsys}) is elliptic whereas the one 
acting on $\Phi$ in 
the second equation of (\ref{DSIsys}) is hyperbolic. 
Davey-Stewartson systems  are of  importance in 
applications since they can be seen as simplifications of the 
Benney-Roskes \cite{BR} and Zakharov-Rubenchik \cite{ZR} systems, 
`universal' models for the description of the interaction of short and 
long waves. These equations have first appeared  in the context of water waves 
\cite{DS,DR,AS,Lan} in particular
in the study of the modulation of plane waves.
DS systems also appear  
in  ferromagnetism \cite{Leb}, plasma 
physics \cite{MRZ}, and nonlinear optics \cite{NM}.      
For more details on DS and its applications the reader is referred to 
\cite{KSDS,KSint} where a comprehensive list of references is given. 
Note that DS I is also interesting from a purely mathematical point of view, since it is a 
nonlinear dispersive partial differential equation, and since it is 
one of the few completely integrable equations in two spatial 
dimensions, see \cite{Fokas,FA}. Local existence results for Cauchy 
problems with small 
initial data were proven in \cite{HS,Chi,HH}, and without a smallness 
assumption in \cite{Hay}. 

The hyperbolic form of the second equation in (\ref{DSIsys}) makes it 
convenient to introduce characteristic coordinates
\begin{equation}
    \xi=x-y,\quad \eta=x+y.
    \label{char}
\end{equation}
In these coordinates DS I (\ref{DSIsys}) takes the form of a non-local 
nonlinear Schr\"odinger (NLS) equation, 
\begin{equation}
    i\Psi_{t}+2(\partial_{\xi}^{2}+\partial_{\eta}^{2})\Psi+
    [(\partial_{\xi}^{-1}\partial_{\eta}+\partial_{\eta}^{-1}\partial_{\xi})
    |\Psi|^{2}]\Psi=0
    \label{DSIchar},
\end{equation}
where we have formally inverted the d'Alembert operator in the second 
equation of (\ref{DSIsys}). In order to do so,  
one has to specify boundary conditions at infinity, a problem analytically  discussed in
\cite{AMS} for the multiscales approach to 
the Kadomtsev-Petviasvili (KP) equation (for a numerical 
implementation see \cite{KSM}). In \cite{FS} it was shown that 
\emph{radiating boundary conditions} allow for stable localised 
traveling waves called \emph{dromions} which appear in the long-time 
behavior of the solutions to certain initial value problems for DS I. Another possibility applied in this context are vanishing boundary conditions 
for $\Phi$ in (\ref{DSIsys}) for $\xi,\eta\to-\infty$ (or 
$\xi,\eta\to\infty$).
We  define the operator $\partial_{\xi}^{-1}$ (as is standard for the KP equation) via its Fourier symbol, 
\begin{equation}
    \partial_{\xi}^{-1}=\mathcal{F}^{-1}_{\xi}\frac{1}{ik_{\xi}}=\frac{1}{2}\mathcal{P}
    \left(\int_{-\infty}^{\xi}-\int_{\xi}^{\infty}\right) 
    \label{symbol},
\end{equation}
where $\mathcal{P}$ denotes the principal value and $\mathcal{F_{\xi}}$ the 
Fourier transform in $\xi$ with $k_{\xi}$ being the dual Fourier 
variable, and likewise for 
$\partial_{\eta}^{-1}$.  Note that a consequence of this 
definition is that for  $f\in L^{1}(\mathbb{R})$, one has 
\begin{equation}
    (\partial_{\xi}^{-1} 
    f(\xi))(+\infty) = -(\partial_{\xi}^{-1} 
    f(\xi))(-\infty)=\frac{1}{2}\int_{-\infty}^{\infty}f(\xi)d\xi
    \label{vanish}.
\end{equation}
These  \emph{trivial 
boundary conditions} will be the only ones studied 
 in this paper.

Numerical studies of DS I solutions have been mainly performed for 
radiating boundary conditions, see \cite{WW,NY1,NY2,BB,MFP}. In this 
paper we will perform a similar study, but for the trivial boundary 
conditions (\ref{symbol}). Since no explicit solitons   are 
known for this case, we first construct localized stationary solutions numerically 
and show that they are also exponentially localized. We study 
the stability of these solutions, which we will also call dromions for 
simplicity.

In \cite{KM,KMS} we have shown how to regularize terms of the type 
(\ref{symbol}) arising 
in the context of D-bar equations with a 
hybrid approach: we subtract a singular term for which the 
Fourier transform can be analytically found. The term is chosen in a 
way that what is left is 
smooth \emph{within finite 
numerical precision}, and that its Fourier transform can 
be numerically computed (we work here with 
\emph{double precision} which is roughly of the order of $10^{-16}$). Note that the terms treated in this way in 
\cite{KM,KMS} are less singular (they lead to cusps in the Fourier 
domain, but are bounded) than the simple poles 
considered here. Therefore, the regularization approach for DS I  is 
more important than for DS II if high accuracy is to be achieved. We show that we can 
reach machine precision in the studied examples. Since we want to 
numerically study blow-up scenarios, an approach of high accuracy as 
presented here is crucial in order to obtain reliable results.

With this approach, we first construct numerically localized 
stationary solutions to DS I and propose\\
\textbf{Main conjecture} (Part I):\\
\emph{The DS I equation has stationary solutions $\Psi(\xi,\eta,t)=
Q_{\omega}(\xi,\eta)e^{i\omega t}$ for $\omega>0$, where $Q_{\omega}$ 
can be chosen to have values in $\mathbb{R}^{+}$. The solutions are exponentially 
localised.}\\
It is unknown whether these solutions are ground states for the energy 
(\ref{energyb}). 

Then we study the time evolution of localised perturbations of these 
stationary solutions and initial data from the Schwartz class 
$\mathcal{S}(\mathbb{R}^{2})$ of 
rapidly decreasing smooth functions with a single hump. We find\\
\textbf{Main conjecture} (Part II):\\
\emph{Initial data $\Psi(\xi,\eta,0)\in \mathcal{S}(\mathbb{R}^{2})$ 
with a single hump lead to one of the following 3 cases:\\
- if $\Psi(\xi,\eta,0)=Q_{\omega}(\xi,\eta)$, the DS I solution is 
stationary;\\
- if the mass $||\Psi(\xi,\eta,0)||_{2}^{2}<m_{Q} :=||Q||^{2}_{2}$ 
($Q:=Q_{1}$), the solution is simply dispersed to infinity;\\
- if the mass $||\Psi(\xi,\eta,0)||_{2}^{2}>m_{Q}$, there is a 
blow-up of the $L^{\infty}$ norm of $\Psi$ at a finite time $t^{*}$ 
such that 
\begin{equation}
	\Psi(\xi,\eta,t) = \frac{Q(X,Y)}{L(t)}+\tilde{\Psi}(\xi,\eta,t)
	\label{conj},
\end{equation}
where $X$, $Y$ are defined in (\ref{dyn}),  where 
$||\tilde{\psi}||_{2}<\infty$ for all times, and where 
\begin{equation}
	L(t)\propto t^{*}-t.
	\label{L}
\end{equation}}
This means that as in the DS II case conjectured in \cite{KS}, the 
blow-up is of the type being unstable for standard NLS equations. 

The paper is organized as follows: in section 2 we give a brief 
overview on DS I equations. In section 3 we 
introduce the hybrid approach to compute the anti-derivatives in 
$\xi$ and $\eta$. Localized stationary solutions to DS I are 
numerically constructed in section 4. In section 5, we present the 
time evolution approach and test it at the example of the stationary 
solution. In section 6, we study the time 
evolution of these stationary solutions. General localized initial 
data are studied in section 7. Some concluding remarks are added in 
section 8. 

\section{Basic Facts}
In this section we collect some basic facts on the DS I equation.

We will always study the DS I equation in 
characteristic coordinates, i.e., in the form (\ref{DSIchar}) of a 
non-local NLS equation. Note that the sign of the nonlinearity is not 
important as in the case of the DS II equations, where it 
distinguishes a focusing and defocusing variant of the equation, see 
e.g., \cite{KSint}. For DS I, a change of sign of the nonlinearity 
can be compensated by a change of sign of either $\xi$ or $\eta$ and 
does not affect the behavior of the solutions otherwise. 

The DS I equation is completely integrable and thus has an infinite 
number of formally conserved quantities. In this paper, we will 
consider the $L^{2}$ norm and the \emph{energy}
\begin{equation}
    E=\int_{\mathbb{R}^{2}}^{}d\xi d\eta\left\{|\Psi_{\xi}|^{2}+|\Psi_{\eta}|^{2}
    +\frac{1}{4}\left(\partial_{\eta}^{-1}|\Psi|^{2}\partial_{\xi}|\Psi|^{2}+
    \partial_{\xi}^{-1}|\Psi|^{2}\partial_{\eta}|\Psi|^{2}\right)\right\}
    \label{energyb}.
\end{equation}
This form of the energy has been chosen in accordance to the 
definition of the anti-derivatives (\ref{symbol}). It can be shown by 
direct computation that the energy is conserved in this case. 


The DS I equation is expected to have stationary solutions of the 
form $\Psi(\xi,\eta,t)=Q_{\omega}(\xi,\eta)e^{i\omega t}$, where 
$\omega\in\mathbb{R}^{+}$, and where we get with (\ref{DSIchar}) the 
following equation for $Q$, 
\begin{equation}
    -\omega Q+2(\partial_{\xi}^{2}+\partial_{\eta}^{2})Q+
    [(\partial_{\xi}^{-1}\partial_{\eta}+\partial_{\eta}^{-1}\partial_{\xi})
    |Q|^{2}]Q=0
    \label{Q}.
\end{equation}
 We are interested in localized solutions to this equation. Note that 
 if the solution $Q:=Q_{1}$ of (\ref{Q}) is known for $\omega=1$, the solution 
for arbitrary $\omega>0$ follows from 
$Q_{\omega}=\sqrt{\omega}Q(\sqrt{\omega}\xi,\sqrt{\omega}\eta)$. For the same reasons as for the standard NLS equation, $Q$ can be 
chosen to be real for localized solutions of this equation. Note that there 
is an explicit solution to (\ref{Q}) called dromion \cite{Kono}, 
which reads for $\omega=1$
\begin{equation}
	\tilde{Q} = \frac{1}{4\cosh \xi/2 \cosh\eta/2+e^{(\xi+\eta)/2}}
	\label{Qdrom},
\end{equation}
if radiating boundary conditions at infinity are used, i.e., if 
$$\partial_{\xi}^{-1}\mapsto 
\widetilde{\partial}_{\xi}^{-1}+f(\eta),\quad \partial_{\eta}^{-1}\mapsto \widetilde{\partial}_{\eta}^{-1}+f(\xi),$$
where
\begin{equation}
	f(\xi) = \frac{4}{4(1+e^{\xi})}+\frac{1}{4(1+2e^{\xi})}
	\label{f2}.
\end{equation}
It is remarkable that the dromions are exponentially decaying towards 
infinity in all directions in contrast to the lump solution, 
the localized stationary solution to DS II which has an algebraic decrease 
towards infinity. Furthermore, again in contrast to the 
lump, the dromion is not radially symmetric. Note that it is 
unknown whether there is an exponentially localised solution to 
(\ref{Q}) for trivial boundary conditions at infinity.

Below we present some properties of DS I solutions:
\begin{itemize}
	
\item Translation invariance:  with 
$\Psi(t,\xi,\eta)$ a solution to equation (\ref{DSIchar}), also 
	$\Psi(t+t_{0},\xi+\xi_{0},\eta+\eta_{0})$ is a solution, where $t_{0}$, 
	$\xi_{0}$, and $\eta_{0}$ are real constants.

\item Galilei invariance: with $\Psi(t,\xi,\eta)$ a solution to 
	equation (\ref{DSIchar}), 
	$\Psi(t,\xi-v_{\xi}t,\eta-v_{\eta}t)\exp(\frac{i}{2}(v_{\xi}(\xi-tv_{\xi}/2)+v_{\eta}(\eta-tv_{\eta}/2)))$ 
	with $v_{\xi}$, $v_{\eta}$ real constants is also a solution. 
Thus a stationary localized solution can be seen as a soliton to the 
equation after a Galilei transformation.

\item Scaling invariance: with 
$\Psi(t,\xi,\eta)$ a solution to equation
(\ref{DSIchar}),  $\lambda \Psi(\lambda^{2}t,\lambda 
	\xi,\lambda \eta)$ with $\lambda\in \mathbb{R}/\{0\}$ is also a solution. 
	Note that the $L^{2}$ norm of $\Psi$ is invariant under these 
	rescalings. Therefore NLS equations in 2D with a cubic 
	nonlinearity are called $L^{2}$ critical. It is known that there 
	can be a blow-up in finite time of the $L^{\infty}$ norm of the 
	solution for smooth initial data with sufficiently large $L^{2}$ 
	norm, see \cite{MR04,SS99}. There does not appear to be a theorem 
	on whether DS I or DS II solutions can blow up for generic 
	initial data of sufficient mass. 

\item Pseudo-conformal invariance: with 
$\Psi(t,\xi,\eta)$ a solution to 
	(\ref{DSIchar}), also
	\[
          \frac{1}{t}\Psi(1/t,\xi/t,\eta/t)\exp\left(i\frac{\xi^{2}+\eta^{2}}{t}\right)
        \]
        is a solution. This implies together with the translation 
	invariance of DS I that a stationary localized DS I solution  under a pseudo-conformal transformation becomes a solution with a blow-up in finite time. 
	This has been used in the context of DS II by Ozawa \cite{ozawa} 
	to construct an explicit blow-up solution.  Note that due to the 
	oscillatory terms, the solution will not be in $H^{1}$ after a 
	pseudo-conformal transformation even if the original solution is 
	in $L^{2}$ for all $t$. For standard $L^{2}$ critical 
	NLS equations, this blow-up mechanism is unstable, see \cite{SS99} 
	for references. 
\end{itemize}

The generic blow-up mechanism for NLS solutions is self-similar, 
which means one uses the above scaling invariance in $\lambda$ with a 
time dependent factor $L(t)$ in a \emph{dynamic rescaling},
\begin{equation}
	X = \frac{\xi}{L(t)},\quad Y = \frac{\eta}{L(t)},\quad 
	\tau = \int_{0}^{t}\frac{dt'}{L^{2}(t')}, \quad \psi(X,Y,\tau) = 
	L(t)\Psi(\xi,\eta,t)
	\label{dyn}.
\end{equation}
The dynamically rescaled DS I equation (\ref{DSIchar}) then reads 
\begin{equation}
    i\psi_{\tau}+i\epsilon a(X 
\partial_{X}\psi+Y\partial_{Y}\psi+\psi)+2(\partial_{X}^{2}+\partial_{Y}^{2})\psi+
    [(\partial_{X}^{-1}\partial_{Y}+\partial_{Y}^{-1}\partial_{X})
    |\psi|^{2}]\psi=0
    \label{DSIresc},
\end{equation}
where $a = \partial_{\tau}\ln L$.

In the case of a blow-up, the scaling factor $L(t)$ is chosen in a 
way to keep certain norms constant during the time-evolution, for 
instance the $L^{\infty}$ norm of $\psi$. If the blow-up is reached 
for a finite time $t^{*}$, then $\lim_{t\to t^{*}}L(t)=0$ and 
$\lim_{t\to t^{*}}\tau = \infty$. For $L^{2}$ critical NLS equations, 
it is expected that $\lim_{t\to t^{*}}a(t) = 0$. In this case, 
equation (\ref{DSIresc}) reduces to the equation for the stationary 
solution (\ref{Q}) in the limit which would indicate that the blow-up 
is self-similar with $Q$ giving the blow-up profile. Note, that the 
generic blow-up rate for $L^{2}$ critical NLS is given by (see 
\cite{MR04})
\begin{equation}
    L(t) \propto \sqrt{\frac{t^{*}-t}{\ln|\ln(t^{*}-t)|}}.
    \label{loglog}
\end{equation}
One of the questions to be addressed in this paper numerically is 
whether there is blow-up in DS I solutions, and whether it follows 
the behavior (\ref{loglog}) or the pseudoconformal rate as in DS II, 
see the conjecture in \cite{KS}. To this end we will trace the 
$L^{\infty}$ norm of $\Psi$ and the $L^{2}$ norm of $\Psi_{\xi}$. 
Both are proportional to $1/L(t)$ and can thus be used to identify the 
scaling factor $L(t)$.

\section{Numerical approach for DS I}
\label{Sec:HybNum}
In this section we briefly describe the numerical approach for the DS 
I equation, in particular how the antiderivatives in (\ref{DSIchar}) are 
computed. We will concentrate here on functions in the Schwartz class  $\mathcal{S}$
of smooth, rapidly decreasing functions. 

The Fourier 
transform of a 1D function $f(\xi)$ and its inverse are defined via
\begin{align}
    \hat{f}(k_{\xi}) & = 
	\mathcal{F}_{\xi}f:=\int_{\mathbb{R}}^{}
	e^{-i\xi k_{\xi}}f(\xi)d\xi,
    \label{Fx}
    \\
    f(\xi) & 
	=\mathcal{F}_{\xi}^{-1}f=\frac{1}{2\pi}\int_{\mathbb{R}}^{} 
	e^{i\xi k_{\xi}}\hat{f}(k_{\xi}) dk_{\xi}.
    \nonumber
\end{align}

The 2D Fourier transform of a function $\Phi(\xi,\eta)$
is defined as 
\begin{align}
    \hat{\Phi}(k_{\xi},k_{\eta}) & = 
	\mathcal{F}_{\xi\eta}\Phi:=\int_{\mathbb{R}^{2}}^{}\Phi(\xi,\eta) 
    e^{-i(\xi k_{\xi}+\eta k_{\eta})} d\xi d\eta,
    \label{F}
    \\
    \Phi(\xi,\eta) & =\mathcal{F}^{-1}_{\xi\eta}\Phi=
	\frac{1}{(2\pi)^{2}}\int_{\mathbb{R}^{2}}^{} e^{i(\xi 
	k_{\xi}+\eta k_{\eta})}\hat{\Phi}(k_{\xi},k_{\eta}) dk_{\xi}dk_{\eta}.
    \nonumber
\end{align}

The basic idea of the Fourier spectral method, which we are going to apply 
here, is to express every function in terms of a Fourier series and 
approximate the latter via a truncated Fourier series. This is 
equivalent to approximating the Fourier transform (\ref{Fx}) via a 
\emph{discrete Fourier transform} which can be efficiently computed 
via a \emph{fast Fourier transform} (FFT).  It is well known that the 
Fourier coefficients of an analytic periodic function decrease 
exponentially, and thus the numerical error due to the truncation of 
the series will also decrease exponentially, see for instance the 
discussion in \cite{trefethen}. Thus Fourier spectral methods show  exponential convergence for analytic functions, 
sometimes called \emph{spectral convergence}. Here we 
only consider functions in the Schwartz class which can be 
efficiently treated as smooth periodic functions on sufficiently 
large tori within the chosen finite numerical precision (the function and 
all relevant derivatives have to vanish at the domain boundaries to the chosen
numerical precision, here $10^{-16}$). 

Derivatives of a function $f(\xi)\in \mathcal{S}(\mathbb{R})$, i.e., 
$$f'(\xi)= \mathcal{F}_{\xi}^{-1}(ik_{\xi}\hat{f}(k_{\xi})),$$
can then be approximated as mentioned above by approximating the 
standard Fourier transform via a discrete Fourier transform. However, 
for the antiderivative
$$\partial_{\xi}^{-1}f(\xi)= 
\mathcal{F}_{\xi}^{-1}\left(\frac{1}{ik_{\xi}}\hat{f}(k_{\xi})\right),$$ 
the singular Fourier symbol will not lead to an exponentially 
decreasing numerical error if the Fourier transform is approximated 
via an FFT. Thus we use a 
hybrid approach, a combination of numerical and analytical 
techniques,  similar to the approach in \cite{KMS2} for the DS II equation. Concretely we write 
\begin{equation}
	\mathcal{F}_{\xi}^{-1}\left(\frac{1}{ik_{\xi}}\hat{f}(k_{\xi})\right) 
	=\mathcal{F}_{\xi}^{-1}\left(\frac{\hat{f}(k_{\xi})-\hat{f}(0)\exp(-k_{\xi}^{2}/4)}{ik_{\xi}}\right) +
	\hat{f}(0)\frac{1}{2}\mbox{erf}(\xi),
	\label{dinv}
\end{equation}
where the \emph{error function} $\mbox{erf}(x)$ is defined as 
\begin{equation}
	\mbox{erf}(x)=\frac{2}{\sqrt{\pi}}\int_{0}^{x}\exp(-y^{2})dy
	\label{erf}.
\end{equation}
The error function can be computed to machine precision with the 
techniques of \cite{KShilbert} since the integral is essentially a 
Hilbert transform in Fourier space. But for simplicity we use the 
Matlab implementation of the error function here. 

The first term on the right hand side of (\ref{dinv}) is a smooth 
function if the limit 
\begin{equation}
	\lim_{k_{\xi}\to0}\frac{\hat{f}(k_{\xi})-\hat{f}(0)\exp(-k_{\xi}^{2})}{ik_{\xi}} = \hat{f}'(0)
	\label{limit}
\end{equation}
is taken into account via de l'Hospital's rule. Since 
$\hat{f}'(0)=\int_{\mathbb{T}}^{}i\xi f(\xi)d\xi$, this term can be computed 
again with Fourier techniques (just the sum of $\xi f(\xi)$ sampled 
on the collocation points). In this way the first term on the 
right hand side of (\ref{dinv}) is in the Schwartz class  if $f(\xi)$ is. Thus it can be 
efficiently computed with Fourier techniques on a large enough torus. 
Note that a Gaussian was introduced in 
(\ref{dinv}) in order to have an integrand in the Schwartz class 
to ensure the rapid convergence of the numerical approach. Thus the 
first term in (\ref{dinv}) is computed to machine precision with 
Fourier techniques whereas the second is obtained with a Matlab 
algorithm with the same precision. 

We illustrate the efficiency of the algorithm for some examples: for 
a Gaussian $f(\xi) = \exp(-(\xi+1)^{2})$ we work with $N=2^{9}$ Fourier 
modes on the interval $10[-\pi,\pi]$. In this case the Fourier 
coefficients decrease to machine precision. Note that we have 
considered a shifted Gaussian here in order to have a non-vanishing 
derivative at the origin in (\ref{limit}). The difference between the error 
function (times the factor $\sqrt{\pi}/2$) is of the order of 
$10^{-16}$. If we consider with the same 
numerical parameters $f(\xi)=\sinh(\xi+1)/\cosh(\xi+1)^2$, the 
Fourier coefficients decrease to the order of $10^{-14}$, and the 
difference to the exact antiderivative   $-\mbox{sech}(\xi+1)$ is of the same 
order. 

\begin{figure}[!htb]
\includegraphics[width=0.49\hsize]{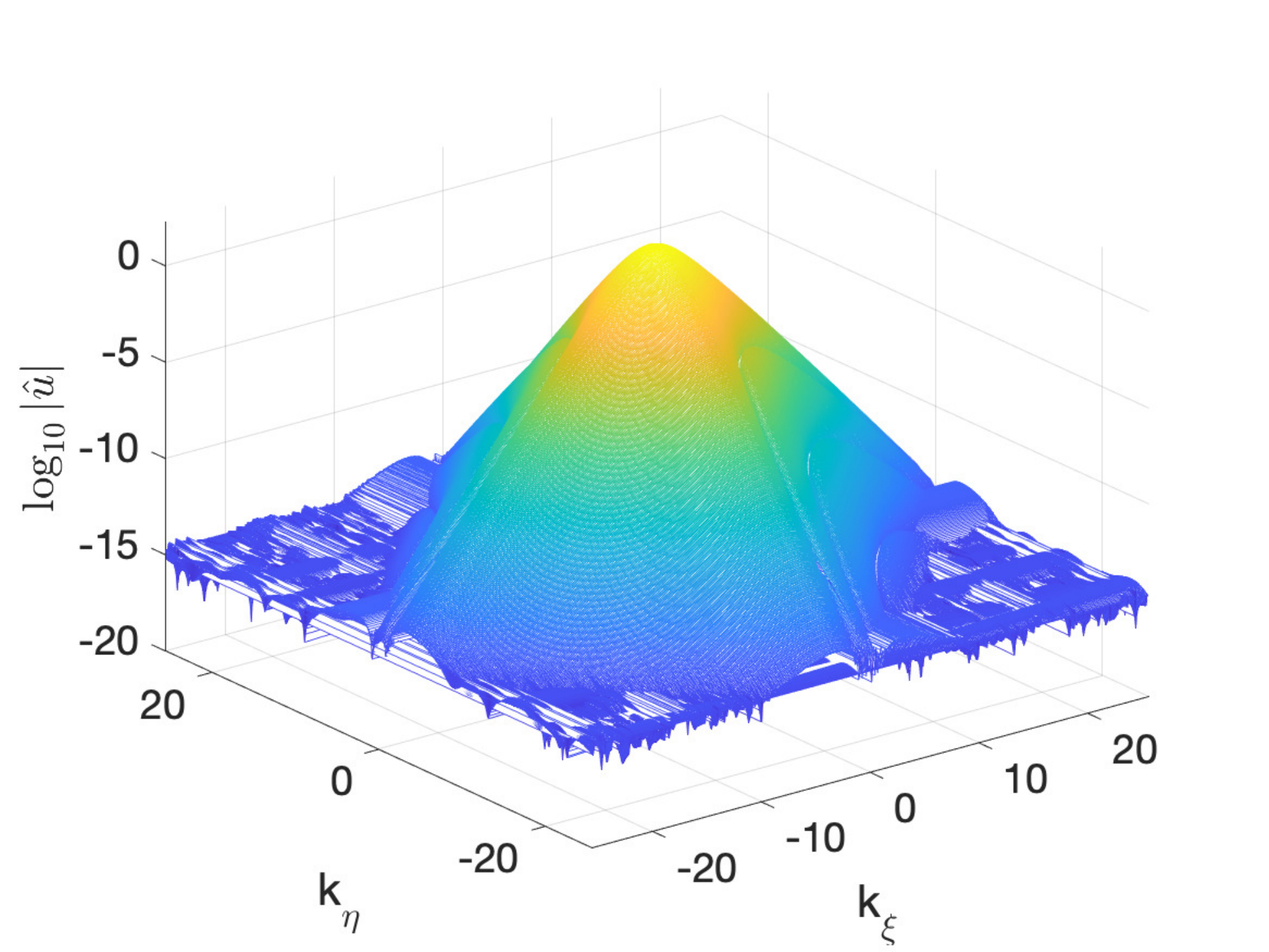} 
\includegraphics[width=0.49\hsize]{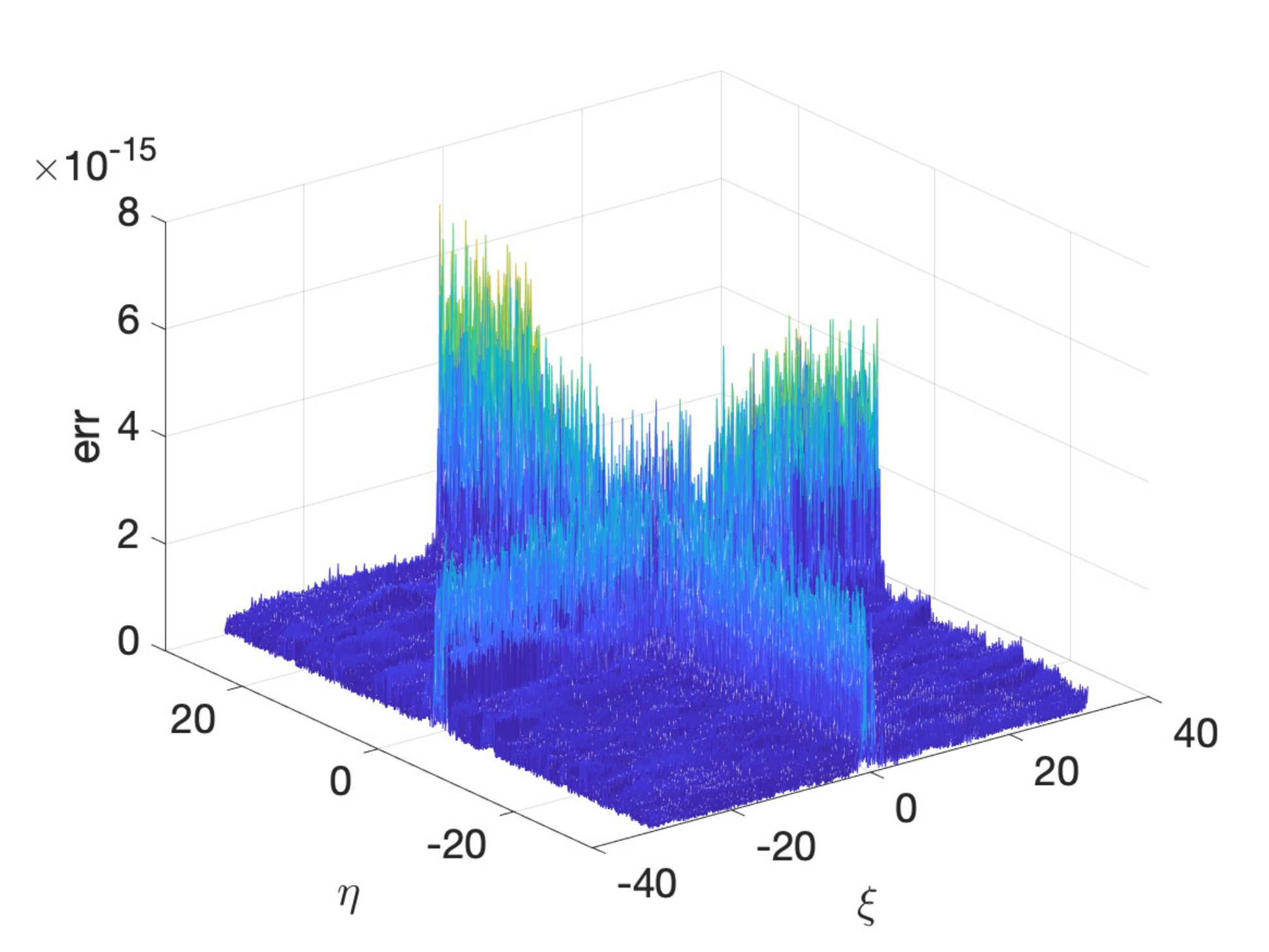} 
\caption{The Fourier coefficients of the function (\ref{dromion2}) on 
the left, and the difference of the numerically computed action of 
the operator 
$\partial_{\xi}^{-1}\partial_{\eta}+\partial_{\eta}^{-1}\partial_{\xi}$ 
on (\ref{dromion2}) and the exact expression on the right. }
\label{dromiontest}
\end{figure}

In the context of DS I, we are obviously mainly interested in the 
accurate numerical computation of the action of the operator 
$\partial_{\xi}^{-1}\partial_{\eta}+\partial_{\eta}^{-1}\partial_{\xi}$ 
on some function in $\mathcal{S}(\mathbb{R}^{2})$. To this end we 
simply apply the above approach in both dimensions. As an example we 
consider the dromion solution $\tilde{Q}_{2}$ for $\omega=2$ in the case of radiating boundary 
conditions,
\begin{equation}
	|\tilde{Q}_{2}|^{2}=\frac{4}{(4\cosh(\xi)\cosh(\eta)+\exp(\xi+\eta))^{2}}
	\label{dromion2}.
\end{equation}
The action of the operator 
$\partial_{\xi}^{-1}\partial_{\eta}+\partial_{\eta}^{-1}\partial_{\xi}$  on the dromion
can be obviously computed explicity. We 
work with $N_{\xi}=N_{\eta}=2^{9}$ Fourier modes in $\xi$ and $\eta$ 
respectively $\partial_{\xi}^{-1}\partial_{\eta}+\partial_{\eta}^{-1}\partial_{\xi}$  on 
$10[-\pi,\pi]\times 10[-\pi,\pi]$. The Fourier coeffcients of 
the function (\ref{dromion2}) can be seen on the left of 
Fig.~\ref{dromiontest}. They decrease to machine precision. The 
difference (denoted by err) between the numerically computed action 
of the operator 
$\partial_{\xi}^{-1}\partial_{\eta}+\partial_{\eta}^{-1}\partial_{\xi}$ on (\ref{dromion2}) 
and the exact expression can be seen on the right of the same figure. 
It is as expected of the same order ($10^{-15}$) as the highest 
Fourier coefficients. 

\section{Localized stationary DS I solutions}
In this section we numerically construct  stationary localized 
solutions to DS I. This is done with the Fourier discretisation introduced in 
the previous section for equation (\ref{Q}) with $\omega=1$. The 
resulting algebraic equation is then iteratively solved with a 
Newton-Krylov method. 

The task is to find a localized solution to (\ref{Q}) where we 
restrict  
ourselves to $\omega=1$ without loss of generality. In Fourier 
space, equation (\ref{Q}) reads
\begin{equation}
    (1 +2k_{\xi}^{2}+2k_{\eta}^{2})\hat{Q}=
    \mathcal{F}_{\xi\eta}\left([(\partial_{\xi}^{-1}\partial_{\eta}+\partial_{\eta}^{-1}\partial_{\xi})
    |Q|^{2}]Q\right)
    \label{Qfourier}.
\end{equation}
As in the previous section, the Fourier transform is approximated via 
a discrete Fourier transform. This implies that (\ref{Qfourier}) leads 
to an $N_{\xi}N_{\eta}$ dimensional nonlinear equation system of the 
form $F(\{\hat{Q}\})=0$ for $\hat{Q}$ (in an abuse of notation, we 
denote the discrete Fourier transform as the standard Fourier 
transform). This system is solved iteratively with a 
Newton method,
\begin{equation}
	\hat{Q}^{(n+1)} = \hat{Q}^{(n)} 
	-\mbox{Jac}(F)^{-1}|_{\hat{Q}^{(n)}}F(\hat{Q}^{(n+1)})
	\label{newton}.
\end{equation}
The action of the Jacobian on $F$ is computed with the Krylov 
subspace method GMRES \cite{gmres}. Note that the Jacobian has a finite dimensional kernel because of the translation invariance of the DS I 
equation. But the iteration converges nonetheless, just the maximum 
of the resulting solution $Q$ will in general not be at the origin. In the plots 
below we have shifted the maximum back to the origin. 

\begin{figure}[!htb]
\includegraphics[width=0.7\hsize]{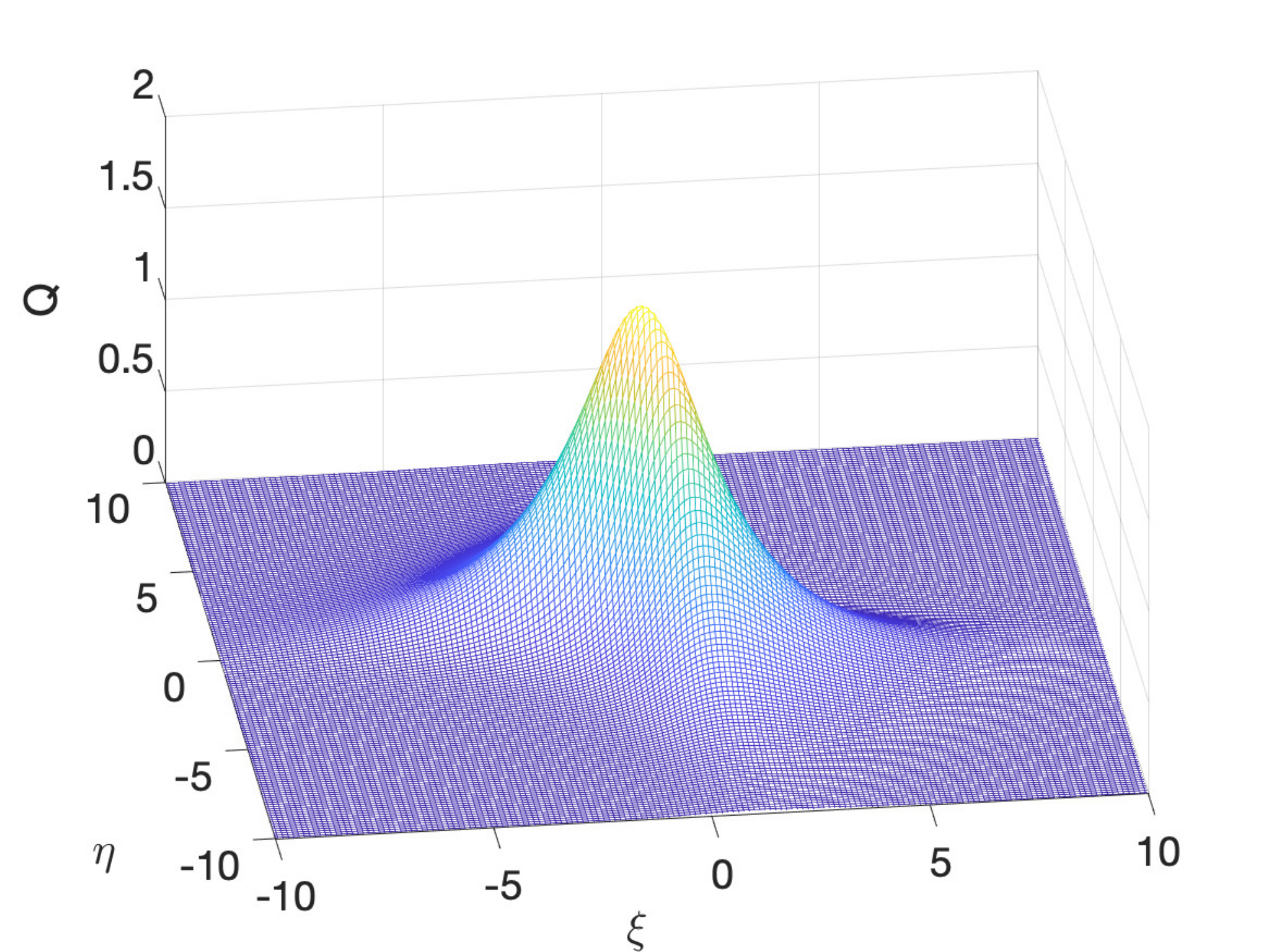} 
\caption{Localized stationary solution to DS I (\ref{Q}) for $\omega=1$. }
\label{dromionn}
\end{figure}

We use $N_{\xi}=N_{\eta}=2^{10}$ Fourier modes for $(\xi,\eta)\in 
20[-\pi,\pi]\times 20[-\pi,\pi]$ and $Q^{(0)} = 
6/(4\cosh(\xi/2)\cosh(\eta/2)+\exp((\xi+\eta)/2))
$ as the initial iterate, i.e., 6 times the dromion (\ref{Qdrom}) for radiating boundary conditions. 
The iteration is stopped once 
$||F||_{\infty}<10^{-10}$. The resulting solution can be seen in 
Fig.~\ref{dromionn}. 

\begin{figure}[!htb]
\includegraphics[width=0.49\hsize]{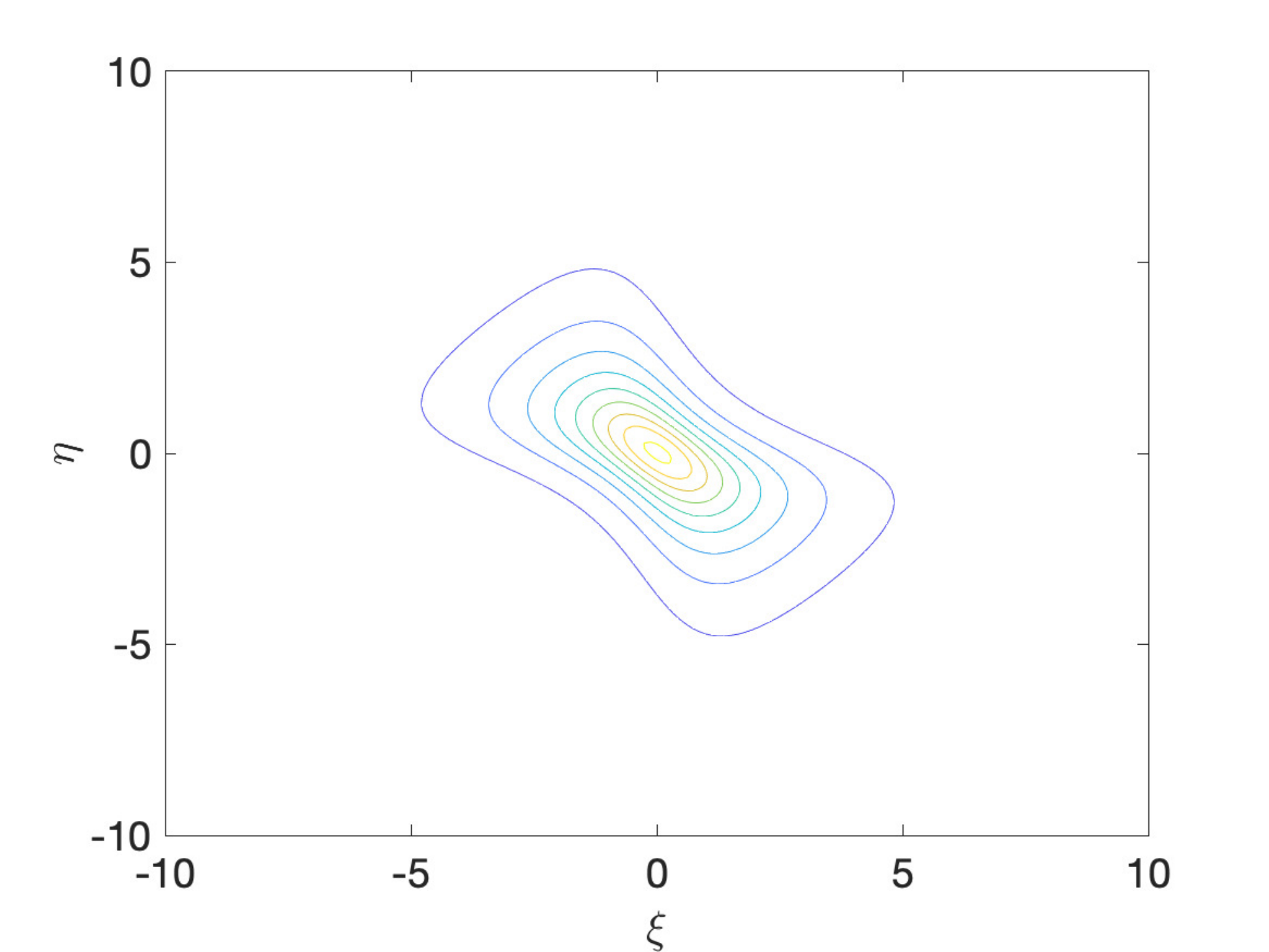} 
\includegraphics[width=0.49\hsize]{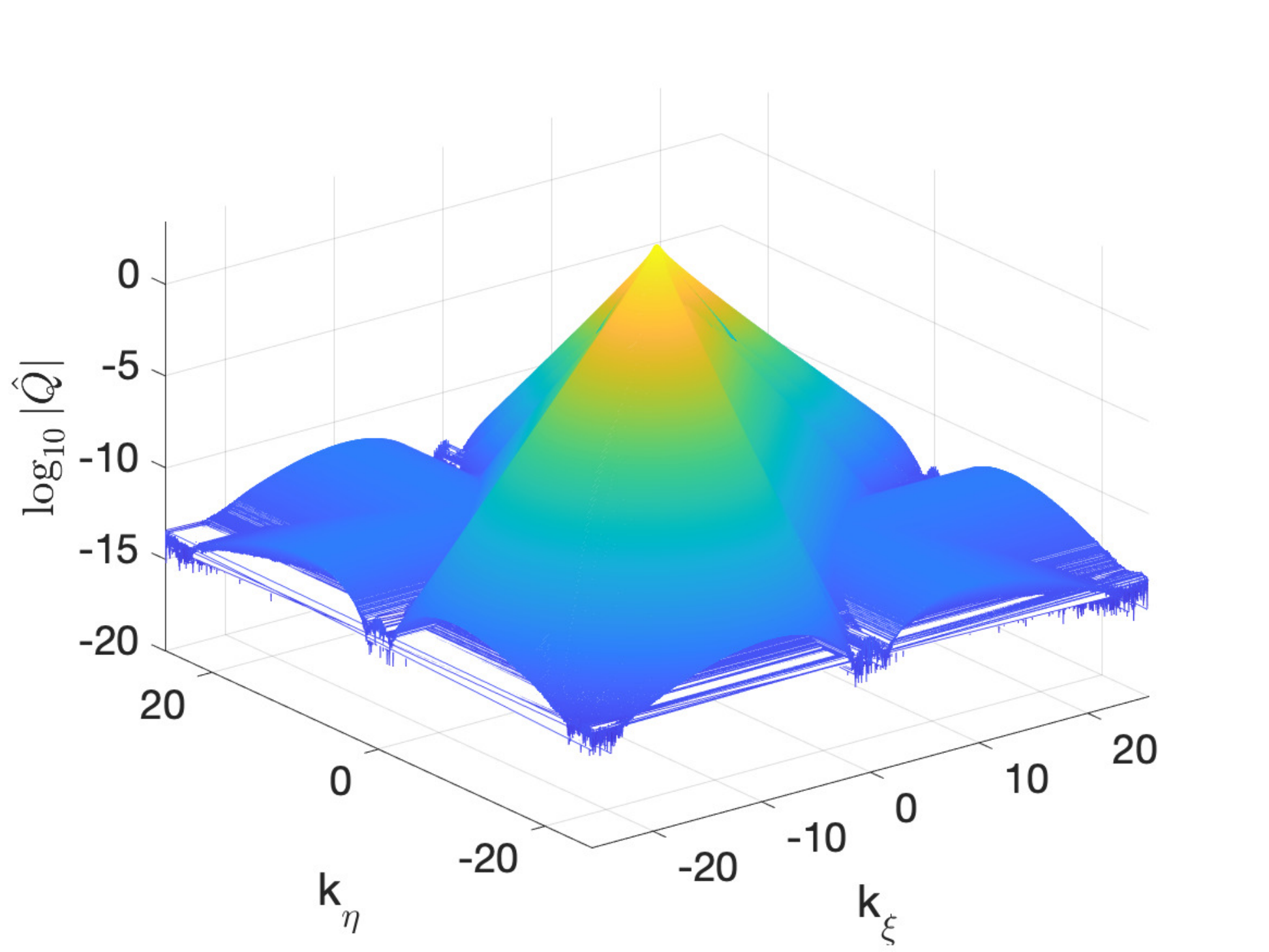} 
\caption{Contour plot of the solution in Fig.~\ref{dromionn} on the 
left, and its Fourier coefficients on the right. }
\label{dromioncontour}
\end{figure}

The solution is again not radially symmetric, but has a symmetry 
with respect to an exchange of $\xi$ and $\eta$ as can be clearly 
seen from the contour plot on the 
left of Fig.~\ref{dromioncontour}. The Fourier coefficents of the 
solution on the right of the same figure decrease to machine 
precision and thus indicate that the solutions is numerically well 
resolved. In fact the numerical parameters have been chosen in a way 
to ensure this. 

\begin{figure}[!htb]
\includegraphics[width=0.32\hsize]{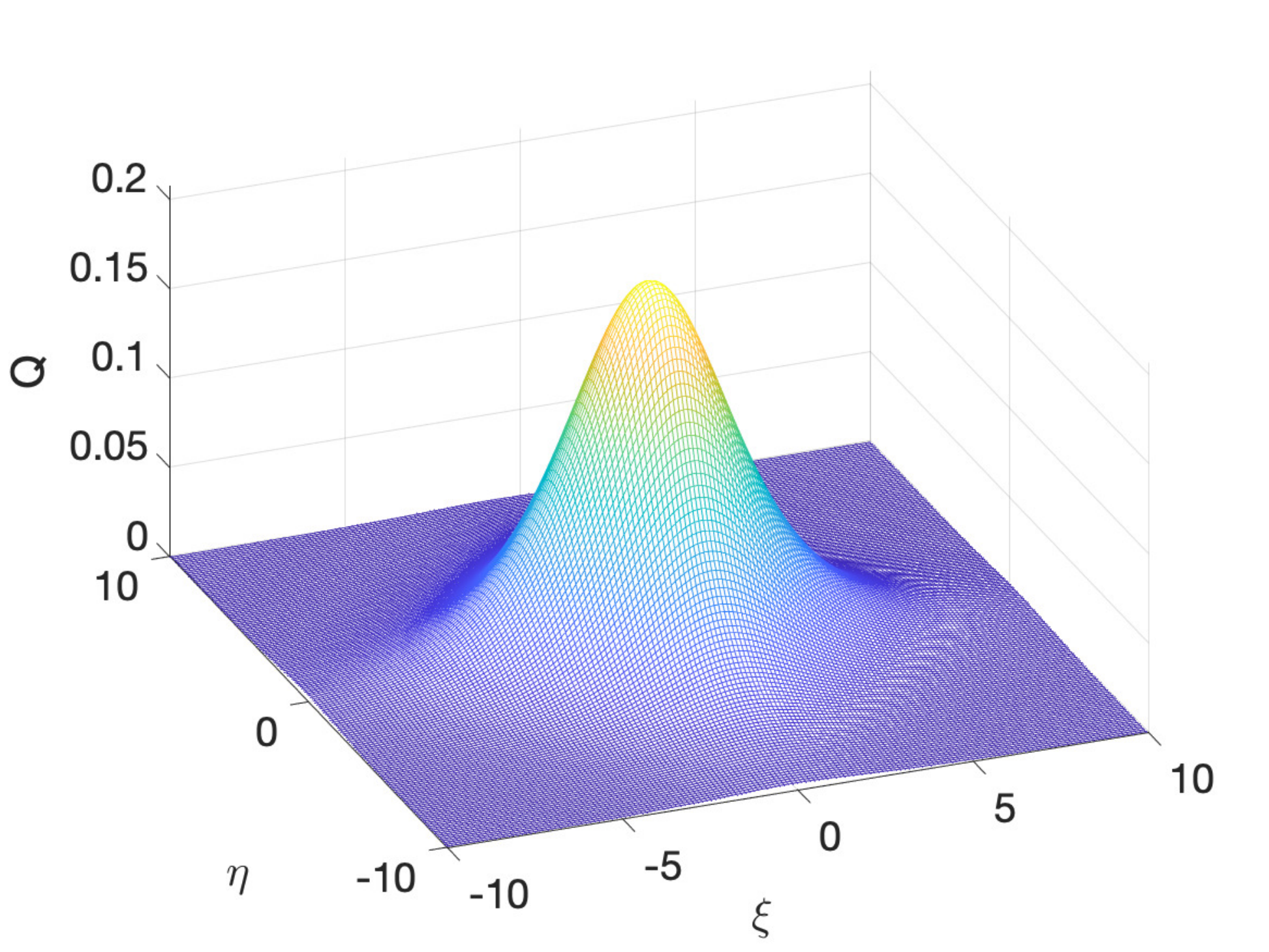} 
\includegraphics[width=0.32\hsize]{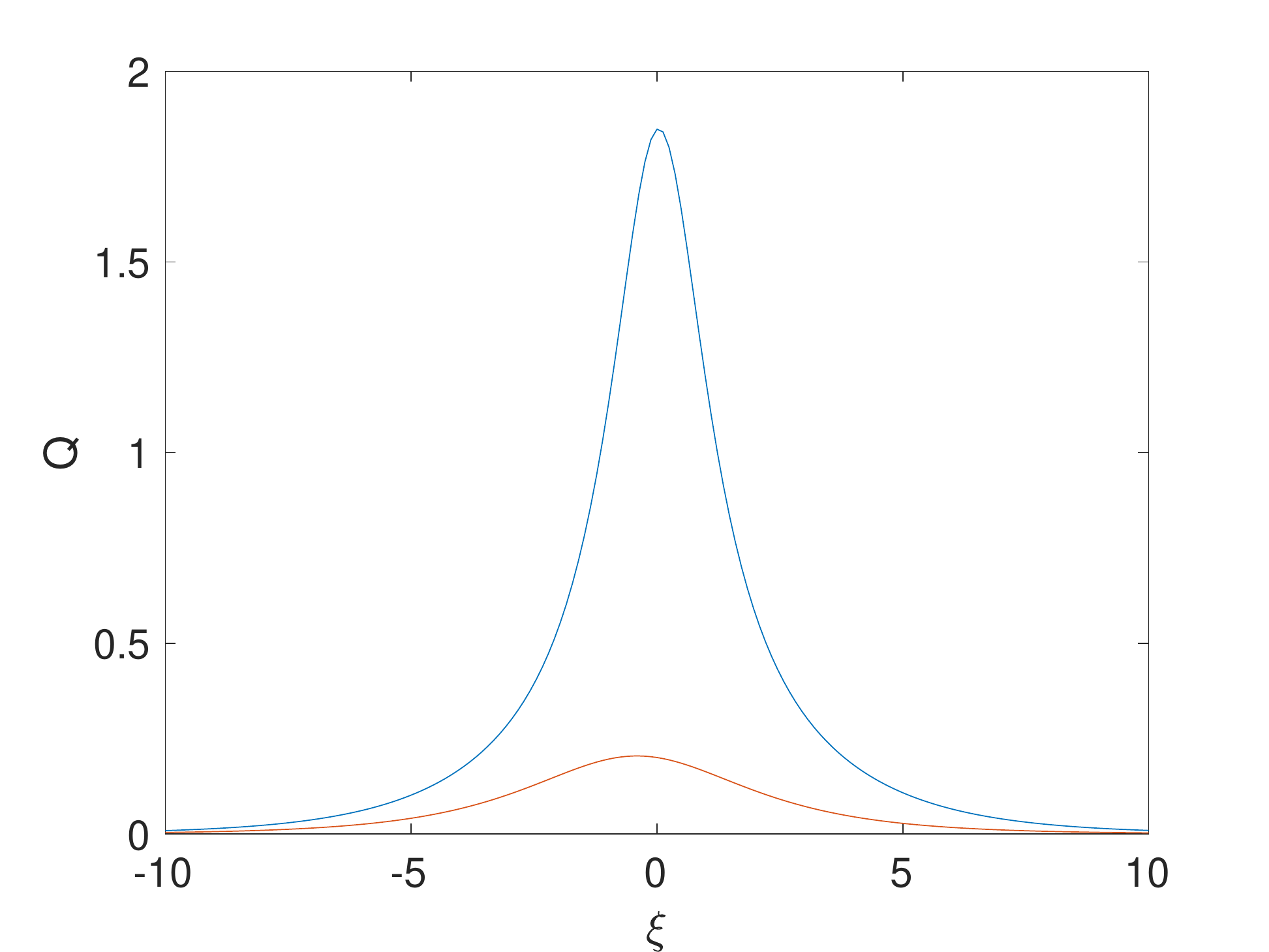} 
\includegraphics[width=0.32\hsize]{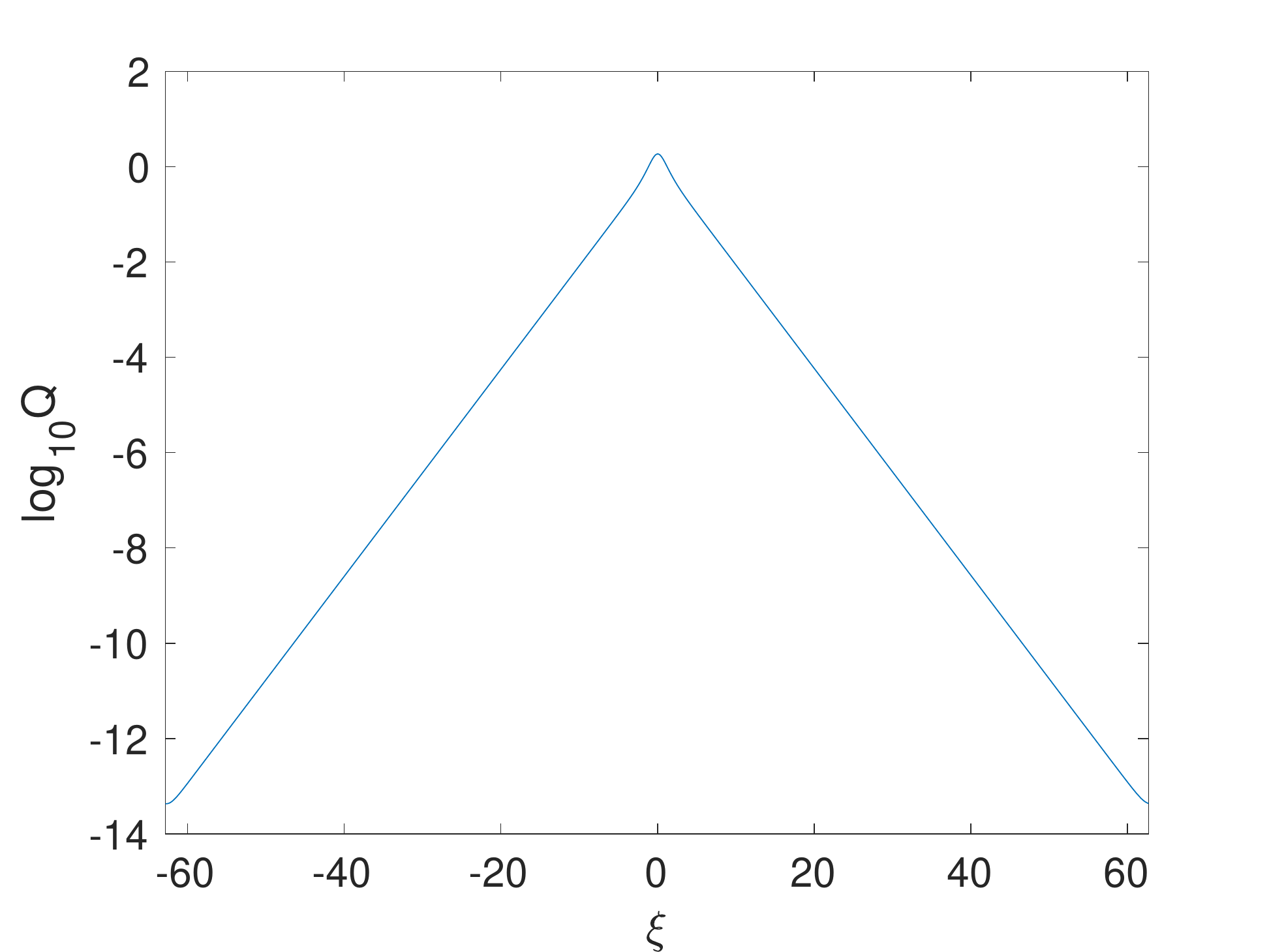} 
\caption{Dromion solution (\ref{Qdrom}) for radiating boundary 
conditions on the left, the same solution on the $\xi$-axis 
(red) together with the 
stationary solution of Fig.~\ref{dromionn} (blue) in the middle, and 
a logarithmic plot of the solution of Fig.~\ref{dromionn} on the 
$\xi$-axis on the right.  }
\label{dromion3}
\end{figure}

Note that the solution is much more peaked than the corresponding one 
(\ref{Qdrom}) for radiating boundary conditions which can be seen on 
the left of Fig.~\ref{dromion3}. In the middle of the same figure, we 
show the solution of Fig.~\ref{dromionn} and (\ref{Qdrom}) on the 
$\xi$-axis in one figure. Obviously the solution constructed in this 
section has a considerably larger maximum (which is why the initial 
iterate had to be chosen with a factor of 6). It is also more slowly 
decaying. However, it is also exponentially decaying as can be seen 
from the logarithmic plot on the right of Fig.~\ref{dromion3} on the 
$\xi$-axis. We only show the plot on the $\xi$-axis here, but the 
same behavior is observed for all values of $\eta$, and for the $\eta$-dependence 
for all values of $\xi$. Thus the stationary solutions to 
DS I are exponentially localized in contrast to the lumps of DS II 
which are algebraically decaying, and this not only for radiating 
boundary conditions. Therefore we will call this solution also 
dromion in the following even though it is not identical to classical one 
in (\ref{Qdrom}). Note that the numerical parameters in this section 
have been chosen in a way that both the solution and its Fourier 
coefficients decrease to machine precision. 

\section{Time evolution}
In this section we outline how the time integration of DS I is 
handled for the discretisation in the spatial coordinates explained 
in the previous sections. We discuss how the accuracy of the time 
integration is controlled and test the code for the example of the 
stationary solution of the previous section. 

We work on $R\times\mathbb{T}_\xi\times\mathbb{T}_\eta$ where 
$\mathbb{T}_\xi = \mathbb{R}/(2\pi L_\xi\mathbb{Z})$, $\mathbb{T}_\eta = 
\mathbb{R}/(2\pi L_\eta\mathbb{Z})$. After the FFT discretisation in 
$\xi$ and $\eta$ of the previous sections, the DS I equation 
(\ref{DSIchar}) becomes an $N_{\xi}N_{\eta}$ dimensional system of 
ordinary differential equations of the form (in an abuse of notation, 
we denote the $N_{\xi}\times N_{\eta}$ matrix obtained for 
$\Psi(\xi,\eta)$ with the same symbol)
\begin{equation}
	\hat{\Psi}_{t}=\mathcal{L} \hat{\Psi}+ \mathcal{N}(\Psi) 
	\label{DSIdisc},
\end{equation}
where 
\begin{equation}
	\mathcal{L} = -2i(k_{\xi}^{2}+k_{\eta}^{2}),\quad 
	\mathcal{N} = i\mathcal{F}_{\xi\eta}\left([(\partial_{\xi}^{-1}\partial_{\eta}+\partial_{\eta}^{-1}\partial_{\xi})
    |\Psi|^{2}]\Psi\right)
	\label{LN}.
\end{equation}
The linear part proportional to $\mathcal{L}$ is diagonal and 
\emph{stiff} since it is quadratic in $k_{\xi}$ and $k_{\eta}$, which 
means that explicit time integration schemes are not efficient. For such 
cases there are many efficient time integration schemes, see for 
instance the references in \cite{etna,KR}. Since it was found in 
\cite{KR} that 
Driscoll's composite Runge-Kutta method \cite{Dri} is very 
efficient for DS equations, we apply it also here. 

We use the  relative conservation of the mass to control the accuracy 
in the time integration. Because of unavoidable numerical errors, the 
numerically computed mass will depend on time even though it is a 
conserved quantity. Thus $\Delta = \log_{10}|1 - 
m/m_0|$, where $m_0$ is the initial mass and $m$ the computed mass can be used to 
control the accuracy of the temporal discretisation. Generally 
$\Delta$ overestimates the temporal resolution by two orders of 
magnitude. The relative mass conservation stays well below $10^{-12}$ throughout most of the runs 
and sharply increases close to the time $t^{*}$ of a potential finite 
time blow-up. Such a jump indicates a loss of precision, and we 
generally discard results with a value of $\Delta$ greater than $-3$. 
Note that we could also use the conserved energy (\ref{energyb}) to 
this end, but the anti-derivatives in (\ref{energyb}) make this 
quantity numerically problematic if resolution in Fourier space is 
lost near a blow-up. The effect is worse for the energy than for the 
DS I solution since in the latter, the anti-derivative
$\partial_{\xi}^{-1}$ is multiplied with $\partial_{\eta}$ which has 
a smoothing effect in the space of Fourier coefficients. Thus the 
energy would underestimate the accuracy near a blow-up which is why 
we use only mass conservation in the following, where such problems 
do not appear. 

\begin{figure}[!htb]
\includegraphics[width=0.49\hsize]{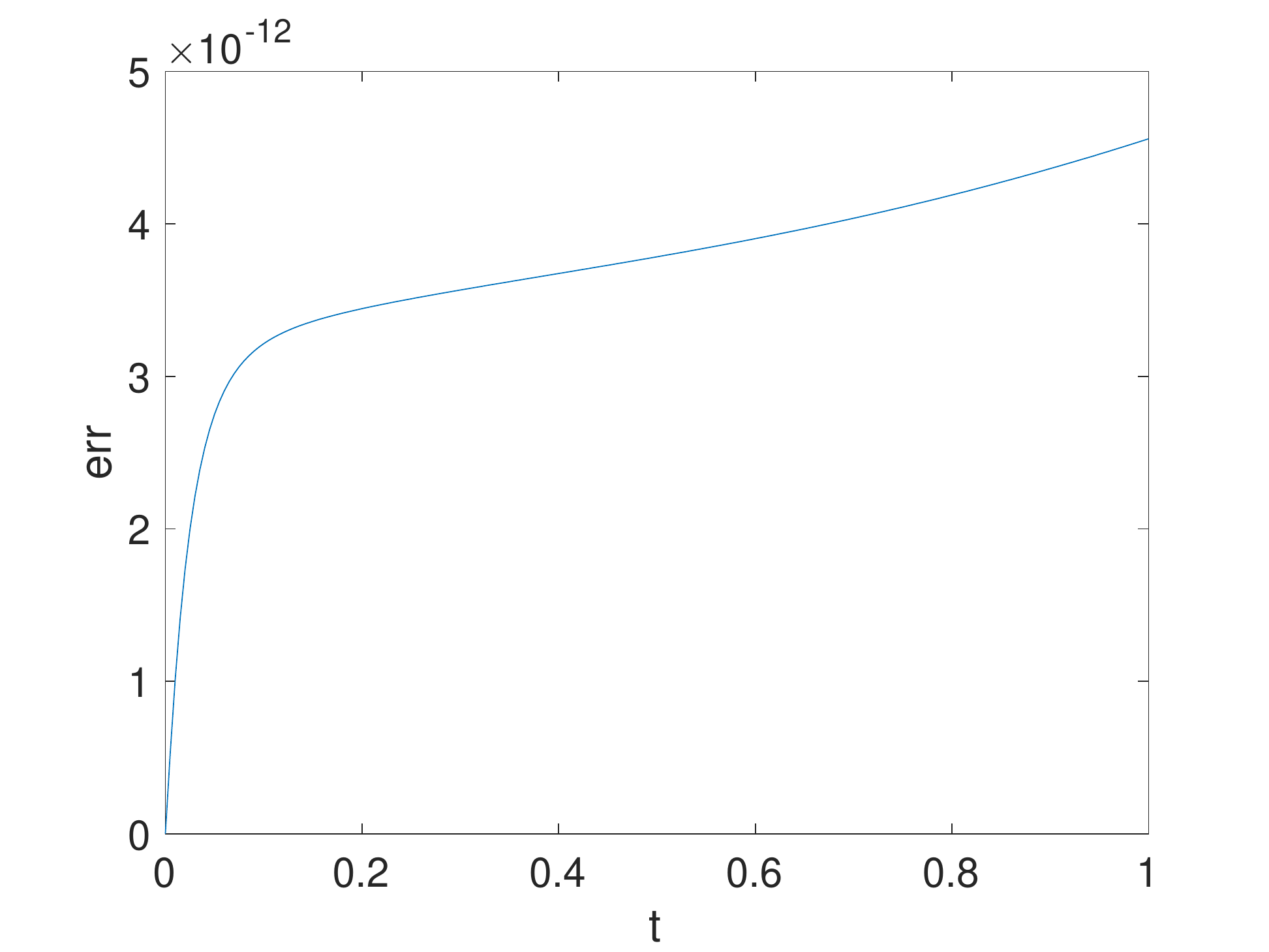} 
\includegraphics[width=0.49\hsize]{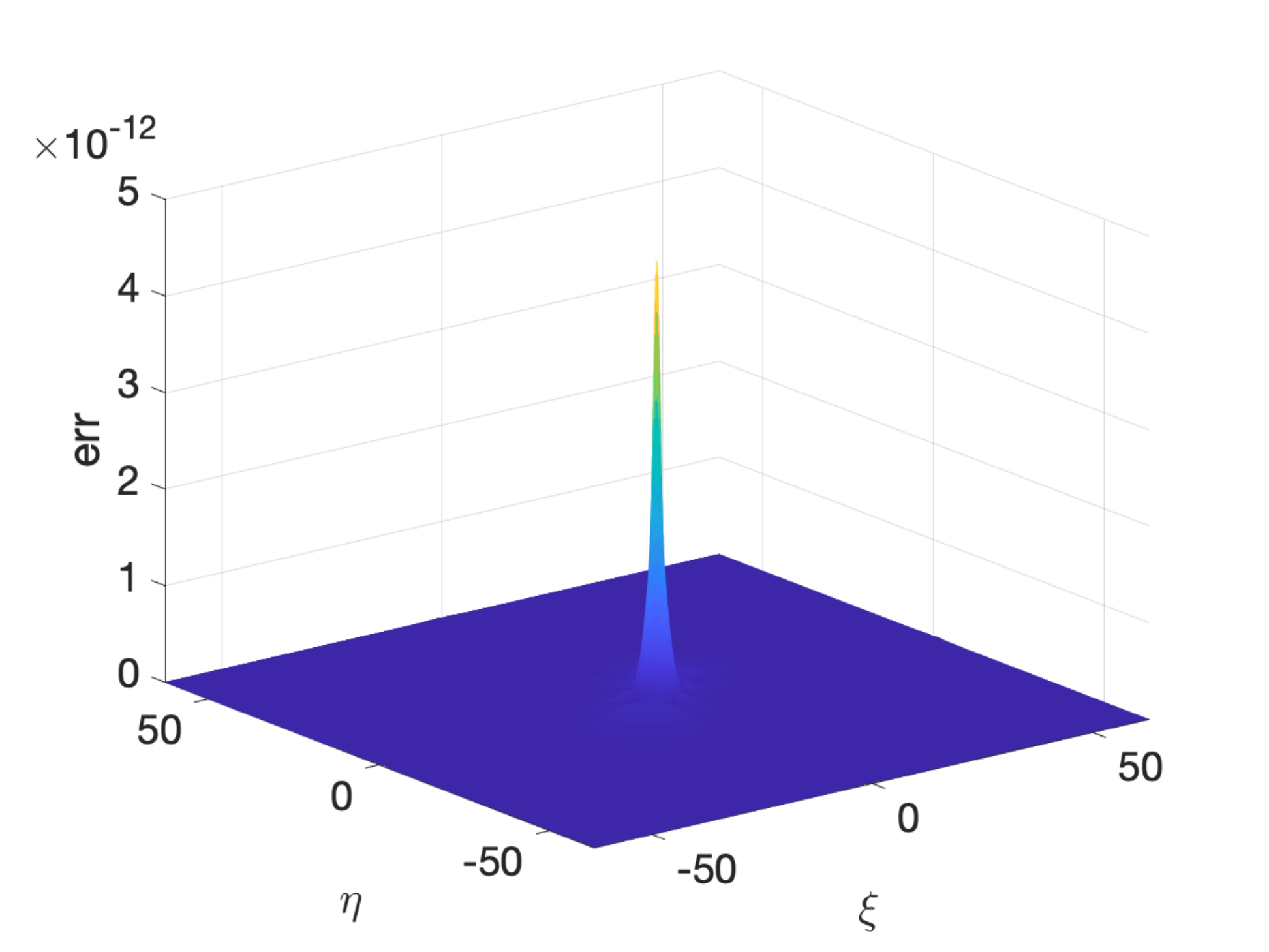} 
\caption{Difference of the DS I solution for the initial data 
$\Psi(\xi,\eta,0)=Q(\xi,\eta)$ and $Qe^{it}$, on the left the 
$L^{\infty}$ norm of the difference in dependence of time, on the 
right the difference for $t=1$. }
\label{dromionerr}
\end{figure}

As an example we consider the dromion constructed in the previous 
section as initial data for DS I. We use $N_{t}=10^{3}$ time steps for $t\leq 1$. The relative 
conservation of the mass is always to the order of $10^{-15}$ (the 
relative energy conservation is of the same order since the solution 
is fully resolved in Fourier space during the whole computation). 
Note that the solution is not static, there is a harmonic time 
dependence. We show the difference between the initial data times 
$\exp(it)$ and the numerical DS I solution in Fig.~\ref{dromionerr}, 
on the left the $L^{\infty}$ norm of the difference between both 
solutions in dependence of time, on the right the modulus of the 
difference for $t=1$ (the difference is always denoted with `err'). 

It can be seen that the difference is during the whole computation of 
the order of $10^{-12}$ or better which is a remarkable result since 
it not only shows the accuracy of the time evolution code, but also 
of the dromion numerically constructed  in the previous section. It 
also shows that the dromion can be stably evolved in time although, as we will show in the following section, it is unstable against perturbations.

\section{Time evolution of the dromion}
In this section we study localized perturbations of the dromion,
mainly of the form 
\begin{equation}
	\Psi(\xi,\eta,0) = \mu Q,\quad \mu>0
	\label{initial}.
\end{equation}
It is shown that perturbations with a mass smaller than the mass of 
the dromion are just dispersed, whereas perturbations with a mass 
larger than the dromion will have a blow-up in finite time. 

\begin{figure}[!htb]
\includegraphics[width=0.49\hsize]{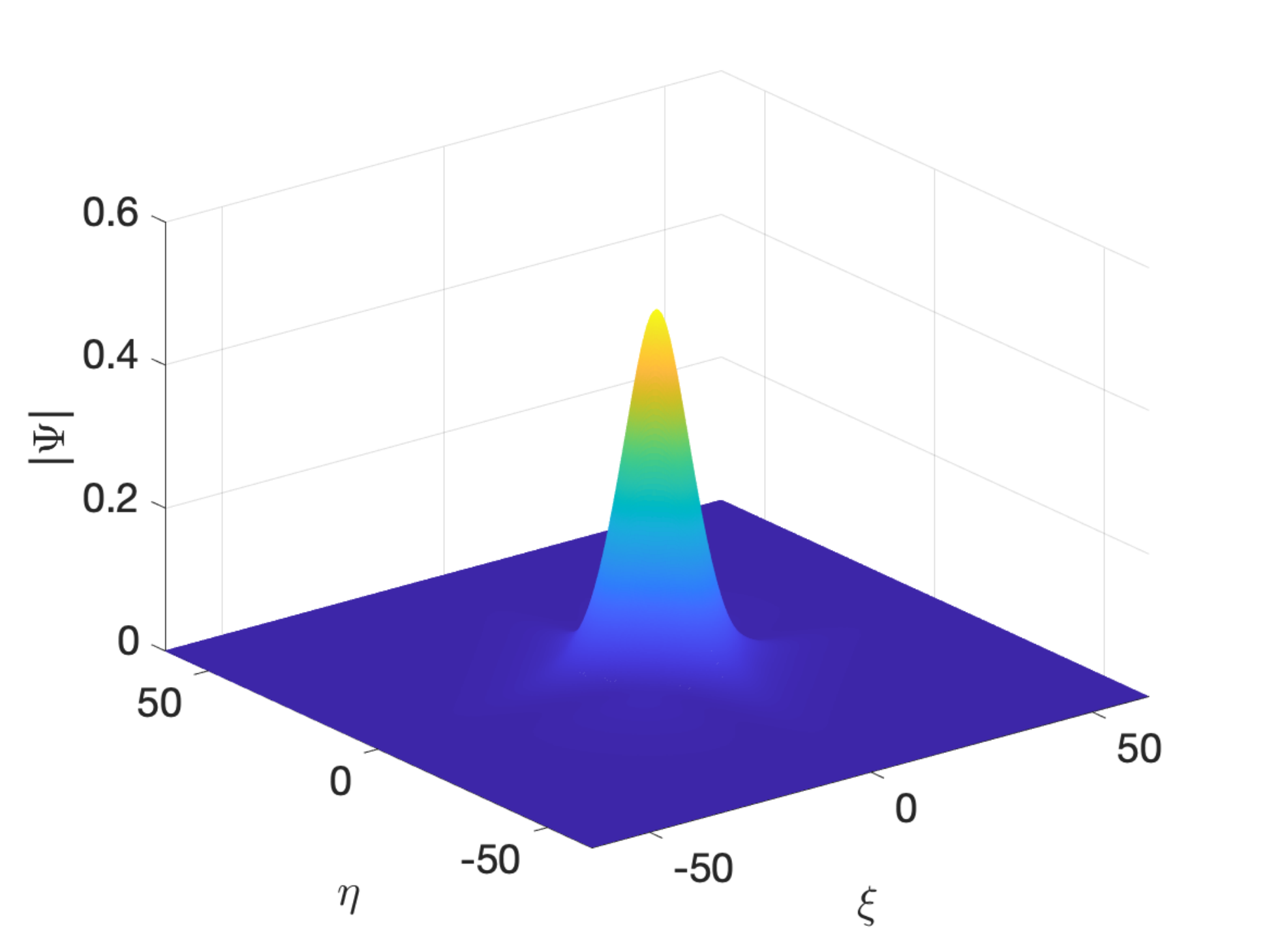} 
\includegraphics[width=0.49\hsize]{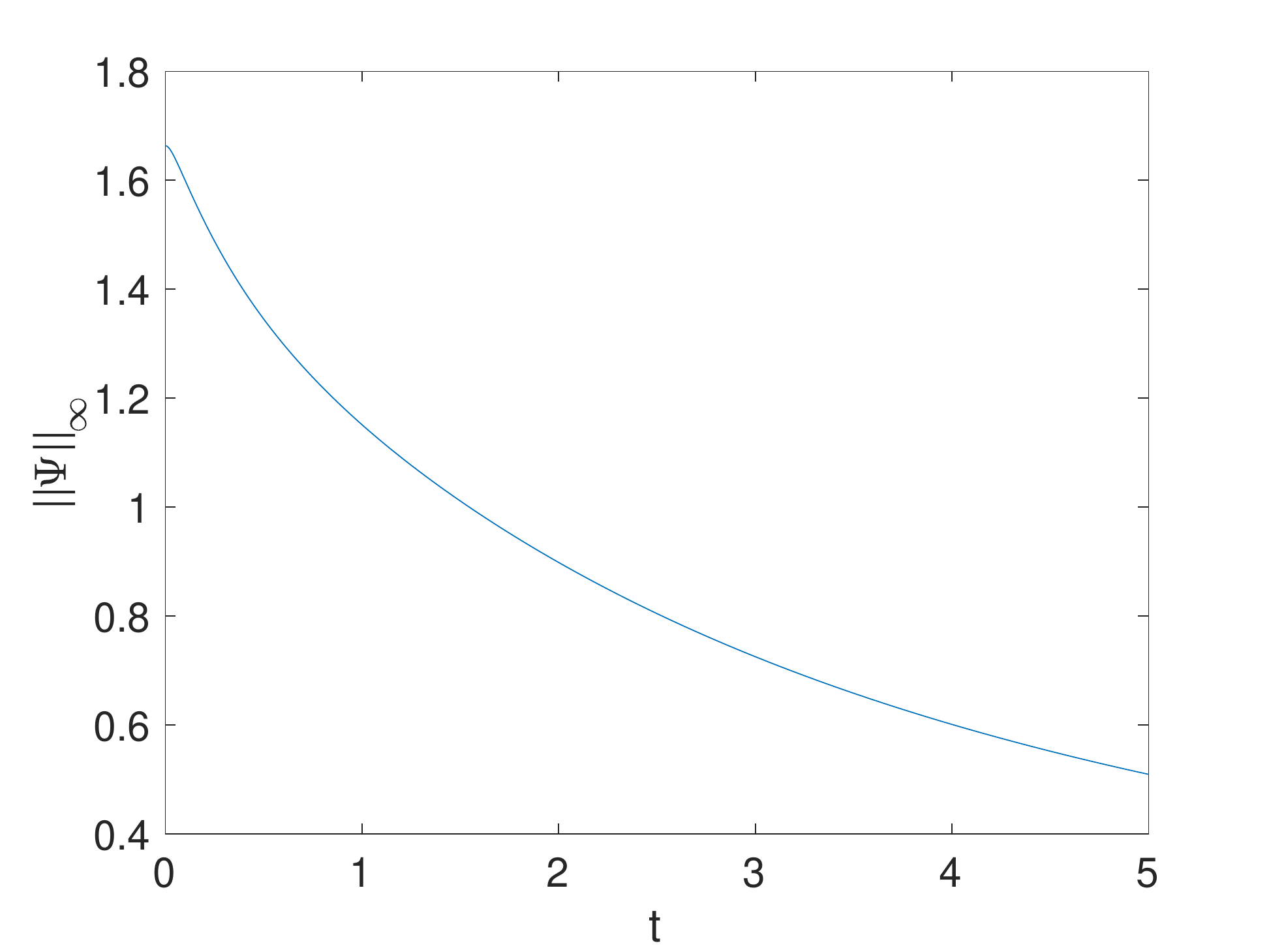} 
\caption{Solution to DS I for the initial data 
$\Psi(\xi,\eta,0)=0.9Q$, on the left for $t=5$, on the right the 
$L^{\infty}$ norm in dependence of time. }
\label{dromion09}
\end{figure}

We first consider the case $\mu=0.9$ in (\ref{initial}) with 
$N_{\xi}=N_{\eta}=2^{10}$ Fourier modes
and $(\xi,\eta)\in 20[-\pi,\pi]\times 20[-\pi,\pi]$, and with 
$N_{t}=5000$ time steps for $t\leq 5$. In this case the initial hump 
simply gets dispersed, it gets wider and flatter over time. The 
solution for $t=5$ can be seen on the left of Fig.~\ref{dromion09}. 
On the right of the same figure the $L^{\infty}$ norm of the solution 
appears to be decreasing monotonically. Note that since we 
approximate a situation on $\mathbb{R}^{2}$ with a setting on 
$\mathbb{T}^{2}$, radiation cannot escape to infinity and thus cannot 
leave the computational domain.  Thus the solution cannot tend to zero 
even for longer times. However we do not find an indication 
of a stable structure in DS I solutions for trivial boundary 
conditions, in contrast to the result in \cite{FS} for radiative 
boundary conditions. 

\begin{figure}[!htb]
\includegraphics[width=0.49\hsize]{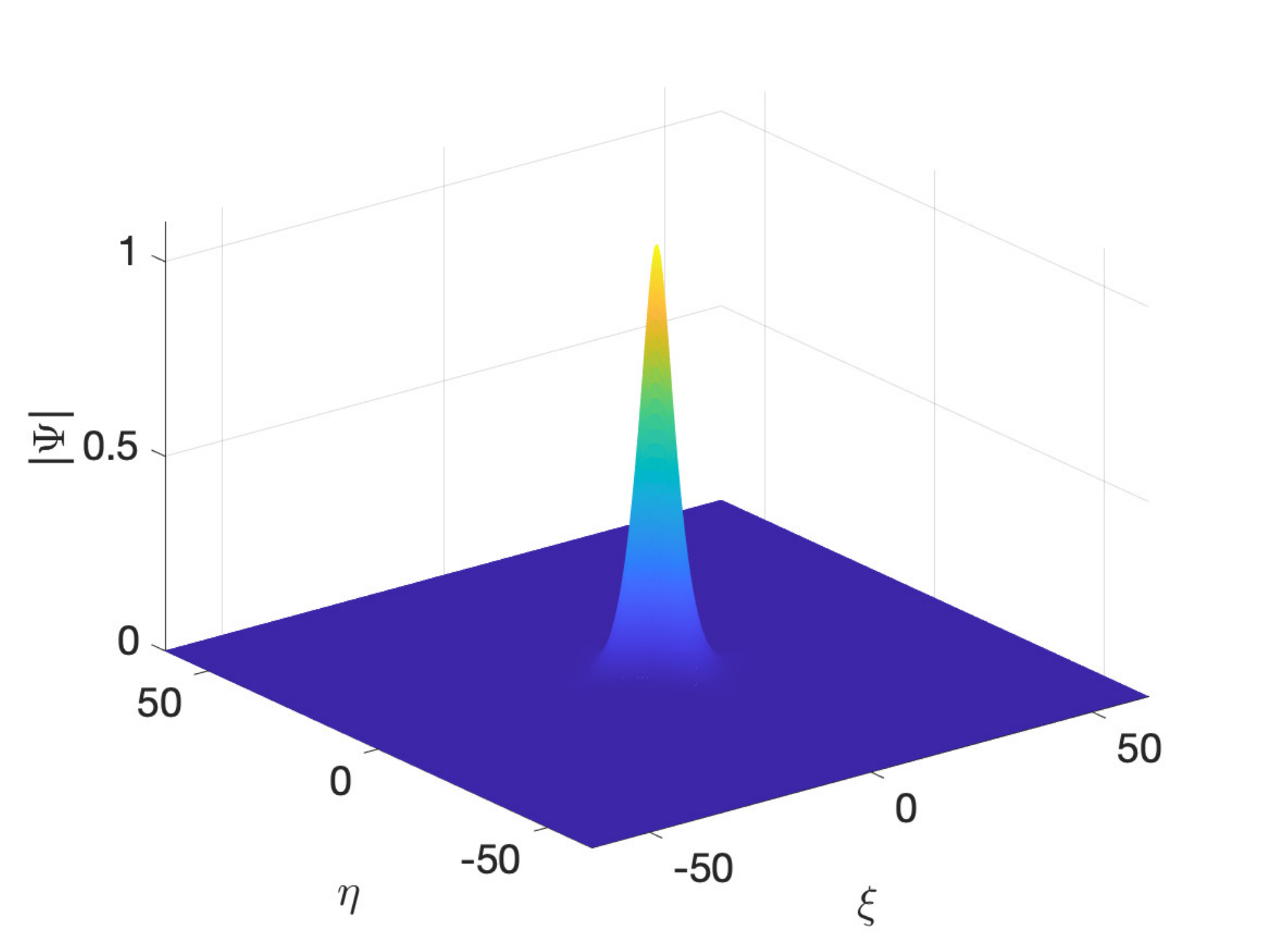} 
\includegraphics[width=0.49\hsize]{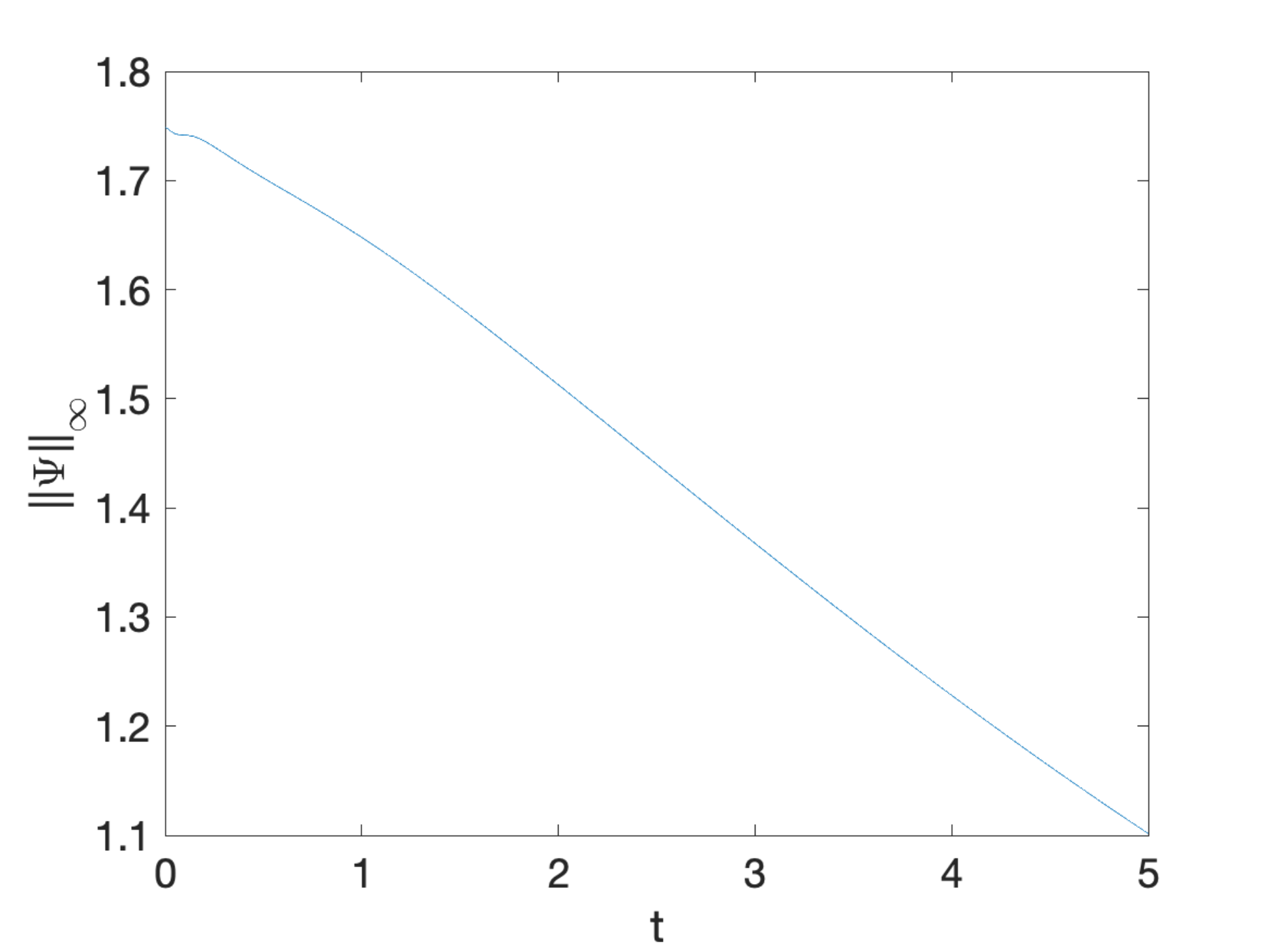} 
\caption{Solution to DS I for the initial data 
$\Psi(\xi,\eta,0)=Q-0.1\exp(-\xi^{2}-\eta^{2})$, on the left for $t=5$, on the right the 
$L^{\infty}$ norm in dependence of time. }
\label{dromionmgauss}
\end{figure}

Then we use the same numerical parameters for the initial data 
	$\Psi(\xi,\eta,0) = Q- 0.1\exp(-\xi^{2}-\eta^{2})$. 
	Note that the mass of these 
data is roughly $0.96 M_{Q}$ where $M_{Q}$ is the mass of the 
dromion and thus larger than the mass in Fig.~\ref{dromion09}. We  
show the solution for $t=5$ on the left of Fig.~\ref{dromionmgauss}. 
Again the initial hump gets just dispersed. This is also confirmed by 
the $L^{\infty}$ norm of the solution on the right of the same figure 
which after some initial oscillation appears to be monotonically 
decreasing. 

The situation changes considerably if we consider perturbations of 
the dromion with larger mass. Here the $L^{\infty}$ norm appears to 
diverge in finite time which obviously cannot be captured 
numerically. However we will trace certain norms in this case and fit 
the found results to the self similar model (\ref{dyn}) for blow-up. This allows us to extend data 
from the region, where the numerical error is still controlled to 
essentially the 
full blow-up scenario, i.e., to identify the blow-up time $t^{*}$ as 
well as the blow-up rate. 

Nonetheless the numerical treatment of a 
blow-up is a delicate problem. 
In order to capture the phenomena, we need to make sure to have 
enough numerical resolution to get close enough to the blow-up in 
order to identify the mechanism. As before, we use 
$L_{\xi}=L_\eta=20$, but now with $N_{\xi}=N_{\eta}=2^{12}$  
Fourier collocation points in each direction. High-index Fourier coefficients, which are 
used to estimate the space resolution, stay below machine precision 
throughout the run. 
In time we use two consecutive runs, one up to $\sim0.9t^{*}$, and a 
second one, with a much finer time step that runs beyond the  
$t^*$ as estimated from the first run. 

\begin{figure}[!htb]\label{dromion11inf}
\includegraphics[width=0.49\hsize]{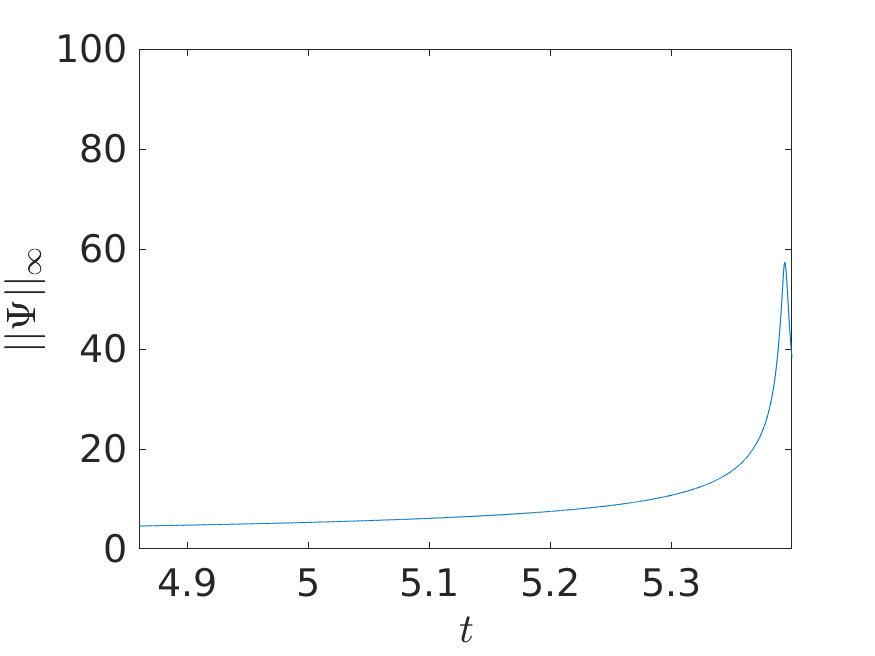} 
  \includegraphics[width=0.49\hsize]{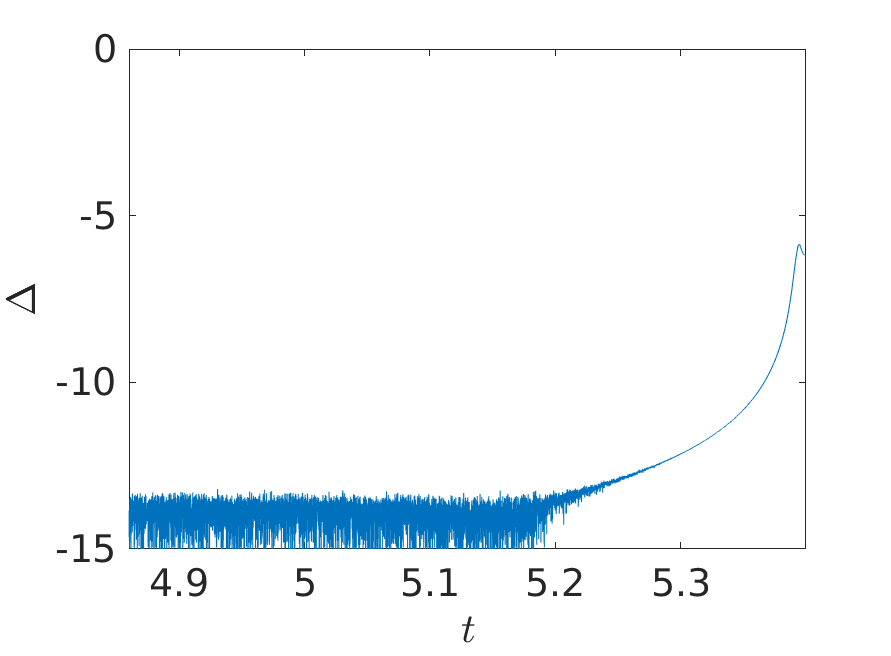} 
\caption{DSI solution close to blow-up for initial data 
$\Psi(\xi,\eta,0) =1.1Q $. On the left the evolution of the 
$L^{\infty}$ norm of the solution, on the right the conservation of 
mass $\Delta = \log_{10}(1-m(t)/m(0))$, which stays below $-14$ until we 
get close to the critical time. The fitted blow-up time is    
$t^* = 5.332$. }
\end{figure}

In Fig.~\ref{dromion11inf} we show on the left  the 
$L^{\infty}$ norm of the solution which appears to indicate a finite 
time blow-up. On the right we trace the quantity $\Delta$ indicating 
the relative conservation of the computed mass. It can be seen to be 
conserved to better than $10^{-5}$ during the whole computation. 

\begin{figure}[!htb]\label{dromion11}
\includegraphics[width=0.7\hsize]{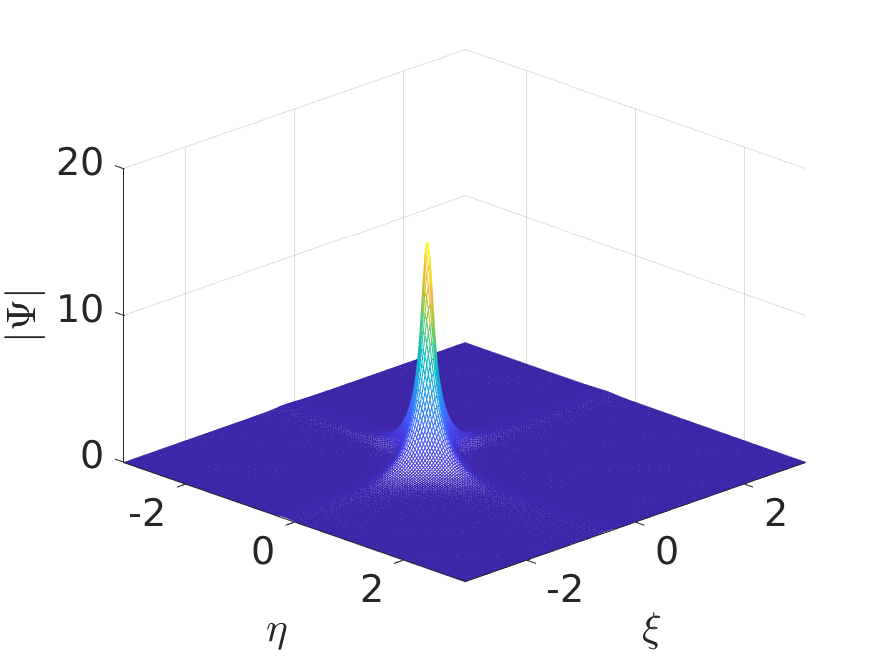} 
\caption{DSI solution close to blow-up for initial data 
$\Psi(\xi,\eta,0) =1.1Q $. }
\end{figure}

We show the solution close to the blow-up in Fig.~\ref{dromion11}. 

\begin{figure}[!htb]
 \includegraphics[width=0.49\hsize]{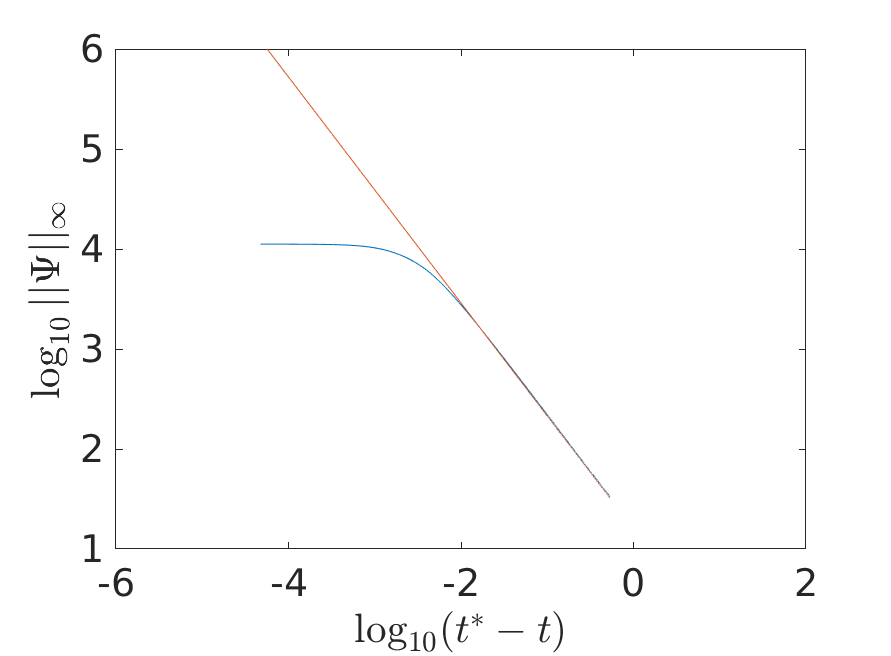}
 \includegraphics[width=0.49\hsize]{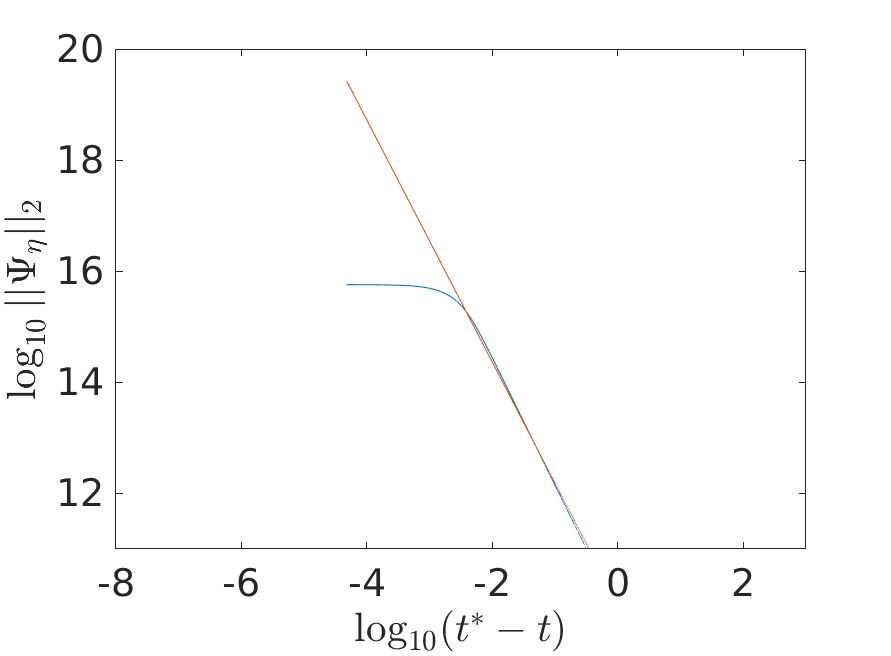}
\caption{Blow-up rate for initial data $\Psi(\xi,\eta,0) = 1.1 Q$. On 
the left $||\Psi||_\infty$, the red lines have the form $a\log_{10}(t*-t) + b$, with values obtained by 
fitting the last several thousand points before we lose precision, $a_\infty = 1.13$, $a_{\Psi_x} = 2.18$ and $b_\infty = 1.2$, $b_{\Psi_x} = 10.0$ and $t^* = 0.5394$. }
\label{Fig:DromBlow}
\end{figure}

To study the mechanism of the blow-up, we trace the $L^{\infty}$ norm 
of the solution as well as the $L^{2}$ norm of the $\xi$ derivative. 
The fitting of these norms is done for the last several thousand 
recorded time steps before we start losing temporal resolution. 
Further we use stabilization of the fit to more precisely judge the 
data cut off point.  Concretely we fit the logarithm of the 
considered norms, for instance the $L^{\infty}$ norm according to $\ln ||\Psi||_{\infty}=a \ln (t^{*}-t)+b$ (and similarly 
for $||\Psi_{\xi}||_{2}$). The 
fitting is performed with the algorithm \cite{NM} implemented in 
Matlab as the command \emph{fminsearch}. For the example with the 
initial data $1.1Q$, the results  can be seen in 
Fig.~\ref{Fig:DromBlow}.

\begin{figure}[!htb]\label{dromion11fit}
 \includegraphics[width=0.7\hsize]{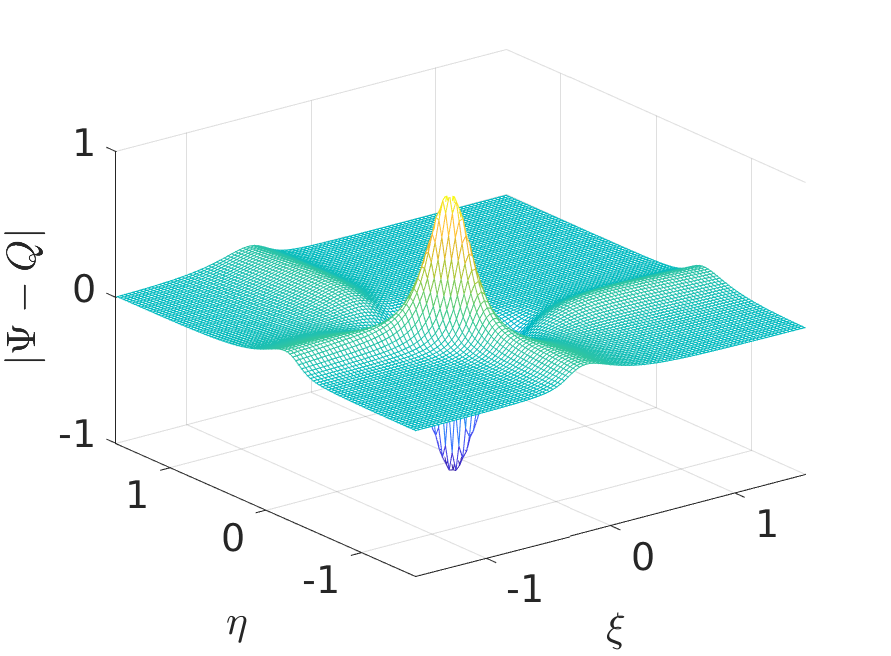} 
\caption{Difference of the DS I solution close to blow-up for initial data 
$\Psi(\xi,\eta,0) =1.1Q $  and a 
scaled dromion.}
\end{figure}

The asymptotic profile of the solution appears to be a scaled dromion 
according to (\ref{dyn}) as can be seen from Fig.~\ref{dromion11fit}. The 
residual is of the order of 10\% of the maximum of the fitted 
solution which shows that one cannot get arbitrarily close to the 
blow-up numerically, but sufficiently to identify the asymptotic 
profile. 


\section{Gaussian initial data}
In this section, we study more examples from the Schwartz class of 
functions with a single hump. Since the dromion is not radially 
symmetric, we concentrate here on standard Gaussians, i.e., initial 
data of the form 
\begin{equation}
	\Psi(\xi,\eta,0) = \kappa \exp(-\xi^{2}-\eta^{2}),\quad \kappa>0.
	\label{inigauss}
\end{equation}
Once more we find that initial data with a mass smaller than the 
dromion will be simply dispersed, whereas initial data with a larger 
mass will lead to a blow-up in finite time. 

\begin{figure}[!htb]
\includegraphics[width=0.49\hsize]{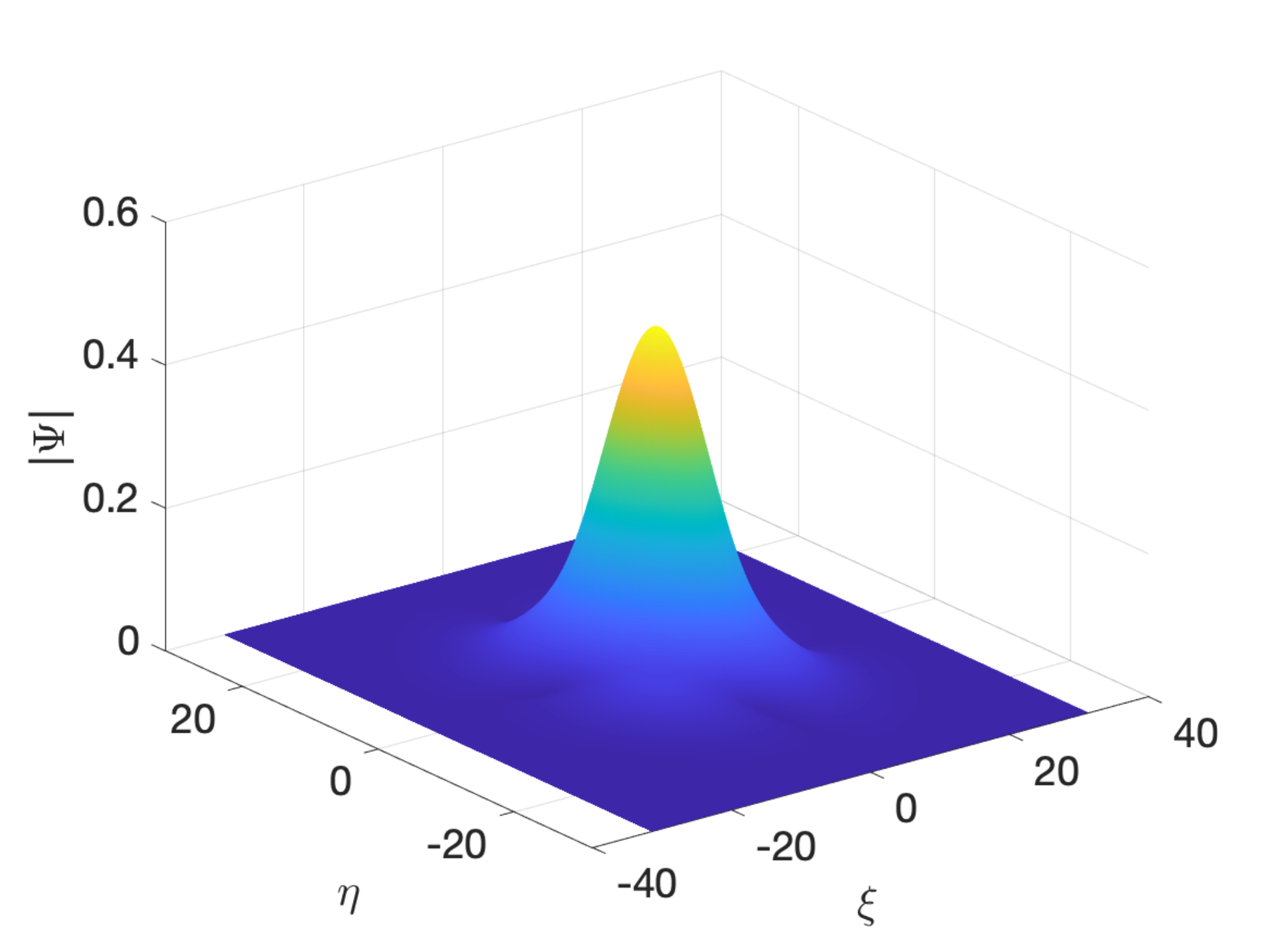} 
\includegraphics[width=0.49\hsize]{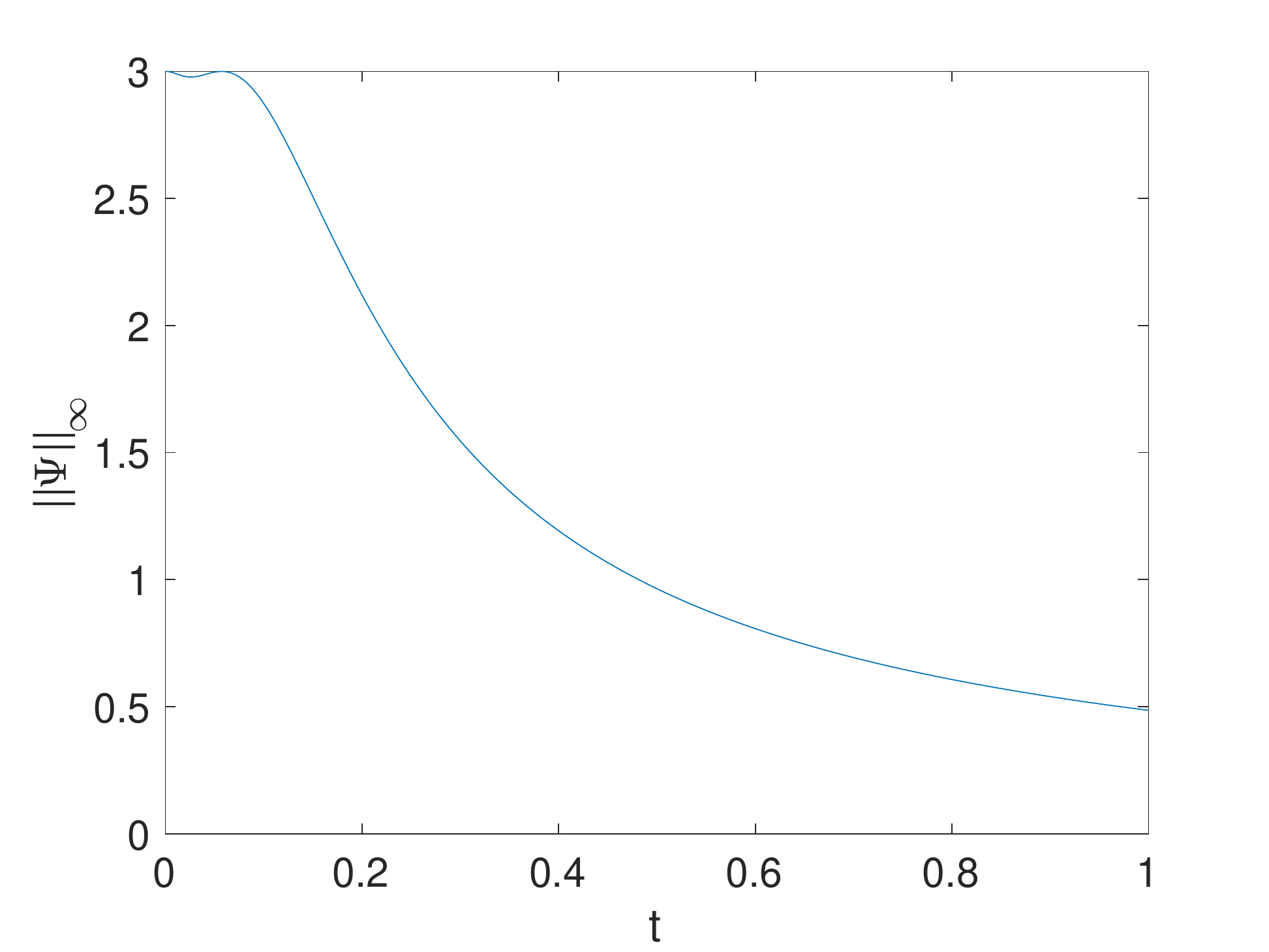} 
\caption{Solution to DS I for the initial data 
$\Psi(\xi,\eta,0)=3\exp(-\xi^{2}-\eta^{2})$, on the left for $t=1$, on the right the 
$L^{\infty}$ norm in dependence of time. }
\label{3gauss}
\end{figure}

First we consider the case $\kappa=3$ in (\ref{inigauss}) with a mass 
of roughly $0.65M_{Q}$. We use $N_{\xi}=N_{\eta}=2^{10}$ Fourier 
modes for $(\xi,\eta)\in 10[-\pi,\pi]\times10[-\pi,\pi]$ and 
$N_{t}=10^{3}$ time steps for $t\leq 1$. The solution for $t=1$ is 
shown on the left of Fig.~\ref{3gauss}. The initial hump is clearly 
 dispersed. This is also confirmed by the $L^{\infty}$ norm of 
the solution on the right of the same figure which after some initial 
growing appears to decrease monotonically. 

\begin{figure}[!htb]\label{5gauss}
\includegraphics[width=0.7\hsize]{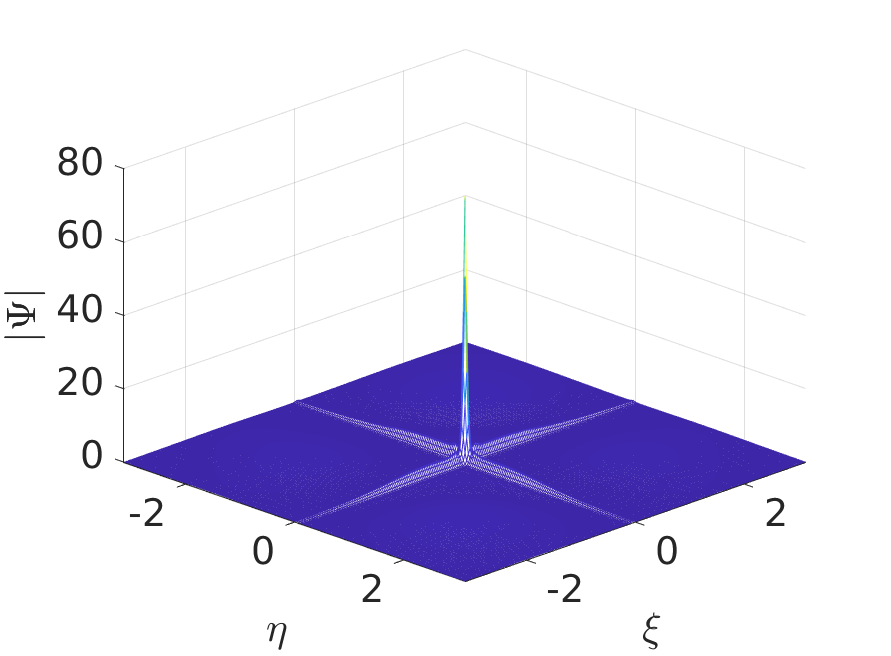} 
\caption{DS I solution close to blow-up for Gaussian initial data 
$\Psi(\xi,\eta,0) = 4.5e^{-\xi^2 -\eta^2}$.}
\end{figure}

If we take initial data with a mass larger than the dromion, say 
$\kappa=4.5$ in (\ref{inigauss}) where the mass is roughly $1.45 M_{Q}$, 
we again seem to get a finite time blow-up. The solution for $t=0.1570$ in 
 Fig.~\ref{5gauss} is already close to the blow-up.

The $L^{\infty}$ norm of the solution is monotonically increasing 
until numerical precision is lost.
\begin{figure}[!htb]
\includegraphics[width=0.46\hsize]{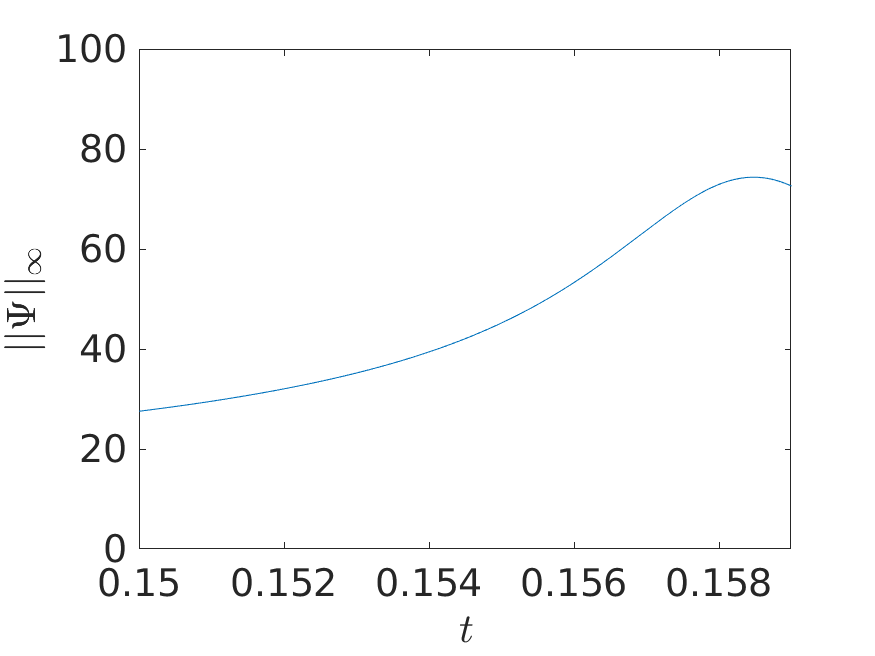} 
 \includegraphics[width=0.49\hsize]{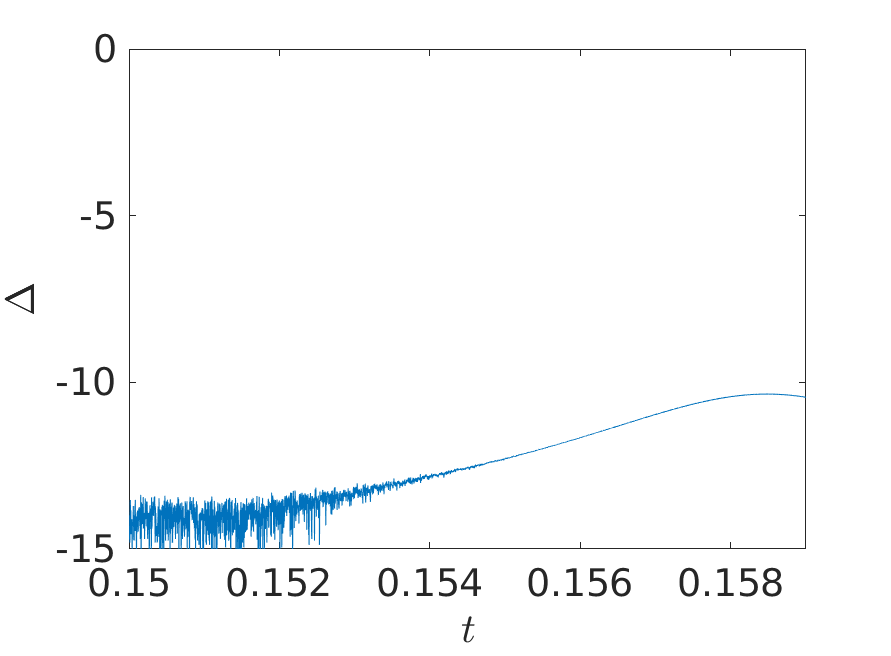} 
\caption{DS I solution close to blow-up for Gaussian initial data 
$\Psi(\xi,\eta,0) = 4.5e^{-\xi^2 -\eta^2}$. Blow up time at $ t = 0.1583$}
\label{Fig:Gauss}
\end{figure}

\begin{figure}[!htb]
 \includegraphics[width=0.49\hsize]{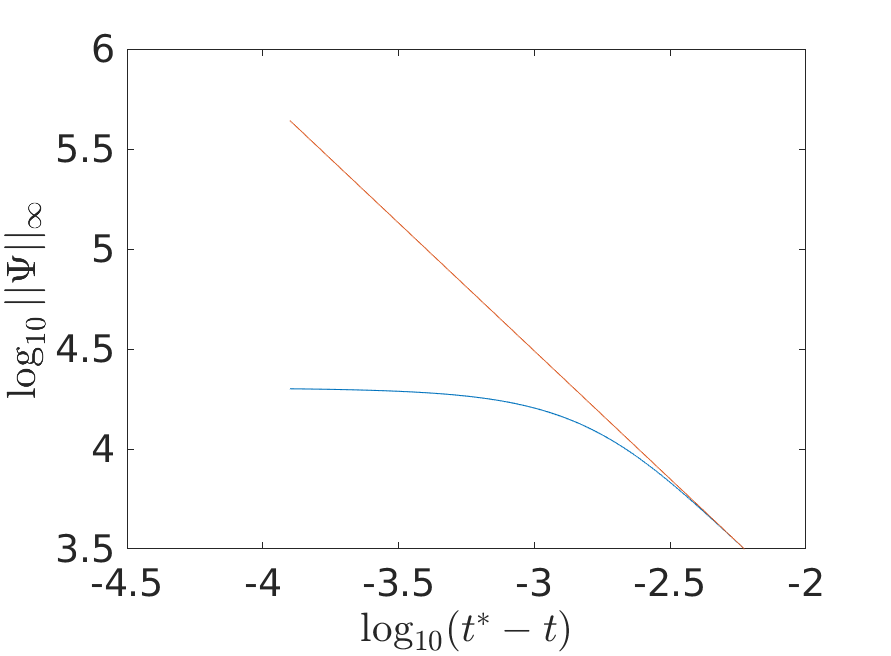}
 \includegraphics[width=0.49\hsize]{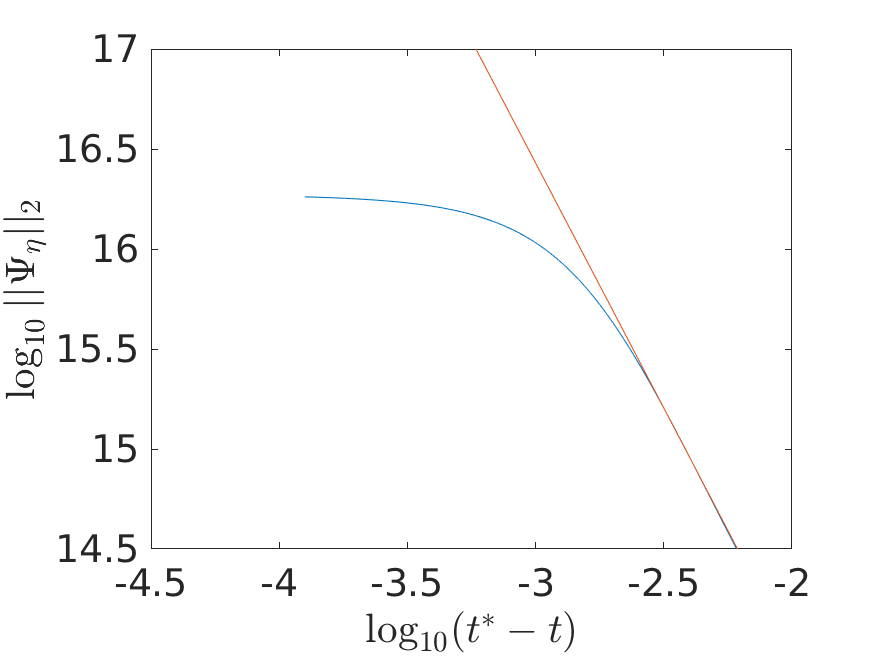}
 \caption{Blow-up rate for initial data $\Psi = 4.5e^{-\xi^2-\eta^2}$, 
 fitting the last several thousand points before we lose precision, $a_\infty = 1.28$, $a_{\Psi_x} = 2.45$ and $b_\infty = 0.65$, $b_{\Psi_x} = 9.08$ and $t* = 0.1583$.}
\label{Fig:GaussBlow}
\end{figure}

A fitting of the $L^{\infty}$ norm of the solution as well as the 
$L^{2}$ norm of $\Psi_{\xi}$ as before in Fig.~\ref{Fig:GaussBlow} 
indicates  that the blow-up is generic, with blow-up rates 
\begin{equation*}
||\Psi||_{\infty} \sim |t^*-t|^{-1}, \quad ||\Psi_{\xi}||_{2}^{2} \approx |t^*-t|^{-2}.  
\end{equation*}   

\begin{figure}[!htb]
 \includegraphics[width=0.7\hsize]{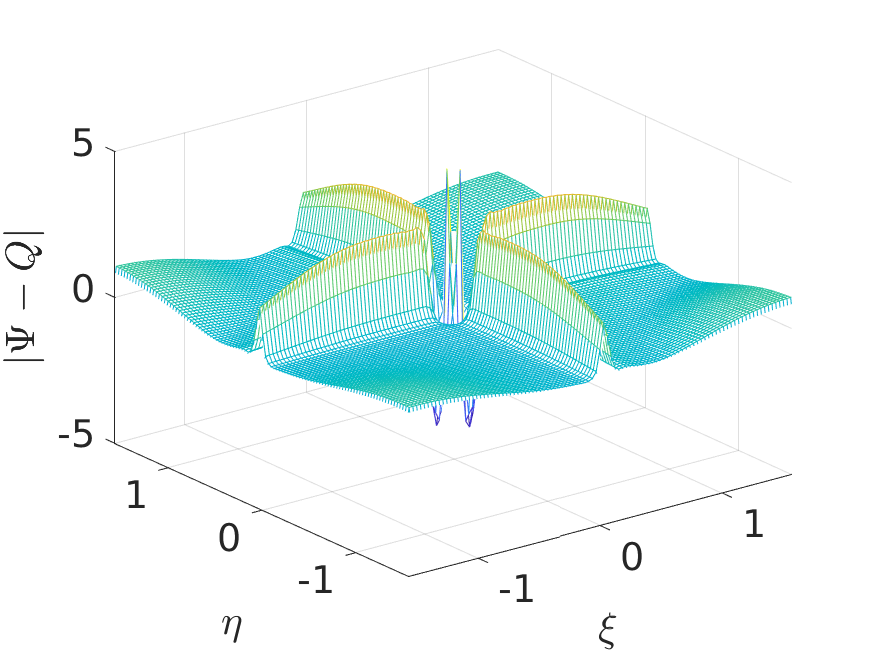}
\caption{Difference of the DS I solution close to blow-up for Gaussian initial data 
$\Psi(\xi,\eta,0) = 4.5e^{-\xi^2 -\eta^2}$ and   a 
fitted (according to (\ref{dyn})) dromion. }
\label{Fig:GaussDiff}
\end{figure}

The blow up is self-similar with the profile close to the blow-up 
being a dynamically rescaled dromion as can be seen from the 
difference between the solution at the final recorded time and 
rescaled dromion (according to (\ref{dyn}) in 
Fig.~\ref{Fig:GaussDiff}. The residual is of the order of $10\%$ 
which once more indicates a good agreement with the model. 

\section{Conclusion}
In this paper, we have presented a detailed numerical study of 
integrable DS I equations with trivial boundary conditions at 
infinity for initial data from the Schwartz class of rapidly 
decreasing smooth functions. As in \cite{KMS2} we have presented a 
hybrid approach based on a Fourier spectral method with an analytic 
(up to the use of the error function) regularisation of the singular 
Fourier symbols. With this approach, it was possible to identify a 
localized stationary solution to DS I which was shown to be 
exponentially localized as the analytically known dromion for 
radiative boundary conditions. Strong numerical evidence has been 
presented that the dromion is 
unstable against localized perturbations, and that perturbations 
leading to a smaller mass of the initial data than the dromion mass 
will be simply dispersed. Perturbations with a larger mass than the 
dromion will lead to blow-up in finite time. We presented numerical 
evidence that the blow-up is self-similar with the dromion as the 
asymptotic profile. The same behavior was observed for initial data 
from the Schwartz class with a single hump. 

An interesting question to be studied in the future is whether 
dromions also exist for non-integrable generalisations of DS I and DS 
II, and 
whether a blow-up is still observed in such cases. A first study of 
these questions for DS II was presented in \cite{KSDS} and should 
also be redone with the methods of \cite{KMS2}.

.
%
%

\end{document}